\documentclass[11pt]{amsart}
\usepackage{mathrsfs}
\usepackage{amssymb}

\usepackage{amscd}
\usepackage{amsmath, amssymb}
\usepackage{amsfonts}
\usepackage[colorlinks,linkcolor=blue,citecolor=blue, pdfstartview=FitH]{hyperref}

\numberwithin{equation}{section}

\theoremstyle{plain}
\newtheorem{thm}{Theorem}[section]
\newtheorem{theorem}[thm]{Theorem}
\newtheorem{lemma}[thm]{Lemma}
\newtheorem{corollary}[thm]{Corollary}
\newtheorem{proposition}[thm]{Proposition}

\theoremstyle{definition}

\newtheorem{question}[thm]{Question}

\newtheorem{remark}[thm]{Remark}

\newtheorem{definition}[thm]{Definition}

\newtheorem{example}[thm]{Example}

\newtheorem{defn-thm}[thm]{Definition-Theorem}

\newcommand{\sE}{{\mathcal E}}

\newcommand{\sO}{{\mathcal O}}

\newcommand{\sX}{{\mathcal X}}

\newcommand{\C}{{\mathbb C}}

\newcommand{\N}{{\mathbb N}}
\renewcommand{\P}{{\mathbb P}}
\newcommand{\Q}{{\mathbb Q}}
\newcommand{\R}{{\mathbb R}}

\newcommand{\T}{{\mathbb T}}

\renewcommand{\S}{{\mathbb S}}

\newcommand{\qtq}[1]{\quad\mbox{#1}\quad}
\newcommand{\bp}{\bar{\partial}}

\newcommand{\ds}{\oplus}

\newcommand{\ts}{\otimes}

\newcommand{\btheorem}{\begin{theorem}}
\newcommand{\etheorem}{\end{theorem}}
\newcommand{\bproposition}{\begin{proposition}}
\newcommand{\eproposition}{\end{proposition}}
\newcommand{\bdefinition}{\begin{definition}}
\newcommand{\edefinition}{\end{definition}}
\newcommand{\bcorollary}{\begin{corollary}}
\newcommand{\ecorollary}{\end{corollary}}
\newcommand{\bproof}{\begin{proof}}
\newcommand{\eproof}{\end{proof}}
\newcommand{\bremark}{\begin{remark}}
\newcommand{\eremark}{\end{remark}}
\newcommand{\eexample}{\end{example}}
\newcommand{\bexample}{\begin{example}}
\newcommand{\la}{\langle}
\newcommand{\elemma}{\end{lemma}}
\newcommand{\blemma}{\begin{lemma}}
\newcommand{\ra}{\rangle}
\newcommand{\sq}{\sqrt{-1}}

\newcommand{\p}{\partial}

\renewcommand{\bar}{\overline}
\newcommand{\eps}{\varepsilon}

\renewcommand{\phi}{\varphi}

\newcommand{\ee}{\end{eqnarray*}}
\newcommand{\be}{\begin{eqnarray*}}

\newcommand{\beq}{\begin{equation}}
\newcommand{\eeq}{\end{equation}}

\newcommand{\bd}{\begin{enumerate}}
\newcommand{\ed}{\end{enumerate}}

\renewcommand{\hat}{\widehat}
\renewcommand{\tilde}{\widetilde}

\renewcommand{\>}{\rightarrow}

\usepackage{fancyhdr}
\renewcommand{\leftmark}{{\scriptsize Big vector bundles and complex manifolds with semi-positive tangent bundles}}

\begin{document}
\title{Big vector bundles and complex manifolds with semi-positive tangent bundles}

\makeatletter
\let\uppercasenonmath\@gobble
\let\MakeUppercase\relax
\let\scshape\relax
\makeatother
\author{Xiaokui Yang}
\date{}
\address{Xiaokui Yang, Department of Mathematics, Northwestern
University, Evanston, IL} \email{xkyang@math.northwestern.edu}

\maketitle

\begin{abstract} We classify compact K\"ahler manifolds with
semi-positive holomorphic bisectional and big tangent bundles.
  We also
classify compact complex surfaces with semi-positive tangent bundles and compact complex $3$-folds of the form $\P(T^*X)$ whose tangent bundles are nef. Moreover, we show that if $X$ is a
Fano manifold such that $\P(T^*X)$ has nef tangent bundle, then $X\cong \P^n$.

\end{abstract}

{\footnotesize{{\textbf{MSC classes}}: $53$C$55$, $32$J$25$,
$32$L$15$, $14$J$15$.}}

\setcounter{tocdepth}{1} \tableofcontents



\section{Introduction}

Since the seminal works of Siu-Yau and Mori   on the solutions to the Frankel conjecture (\cite{Siu-Yau}) and Hartshorne conjecture (\cite{Mo}), it became apparent that positivity properties
of the tangent bundle define rather restricted classes of manifolds. Combining algebraic and transcendental tools, Mok proved the following uniformization theorem in \cite{Mok}: if a compact
K\"ahler manifold $(X,\omega)$ has semi-positive holomorphic bisectional curvature, then its universal cover is isometrically biholomorphic to $(\C^k,\omega_0)\times
(\P^{N_1},\omega_1)\times\cdots \times (\P^{N_\ell},\omega_\ell)\times (M_1, \eta_1)\times\cdots\times (M_k,\eta_k)$ where $\omega_0$ is flat; $\omega_k$, $1\leq k\leq \ell$, is a  K\"ahler
metric on $\P^{N_k}$ with semi-positive holomorphic bisectional curvature; $(M_i,\eta_i)$ are some compact irreducible Hermitian symmetric spaces. Along the line of Mori's work,
Campana-Peternell \cite{CP}  studied the projective manifolds with nef tangent bundles (see also \cite{Yau74} and \cite{Zheng1}); Demailly-Peternell-Schneider (\cite{DPS}) investigated
extensively the structure of compact complex manifolds with nef tangent bundles by using algebraic techniques as well as transcendental methods(e.g. the work \cite{Dem91} of Demailly). For
more details, we refer to \cite{Peternell,CP93,CDP,DPS96} and the references therein.

In the same spirit, Sol\'a Conde and Wi\'sniewski classified projective manifolds with $1$-ample and big tangent bundles:

\btheorem[{\cite[Theoreom~1.1]{SW}}]\label{SW} Let $X$ be a complex projective manifold of dimension $n$. Suppose that the tangent bundle $TX$ is big and $1$-ample. Then X is isomorphic
either to the projective space $\P^n$ or to the smooth quadric $\Q^n$, or if $n=3$ to complete flags $F(1; 2; \C^3)$  in $\C^3$ (which is the same as the projective bundle $\P(T^*\P^2)$  over
$\P^2$).

\etheorem

\noindent A vector bundle $E$ is called \emph{big} if the tautological line bundle $\sO_{\P(E^*)}(1)$ of $\P(E^*)$ is big. The $1$-ampleness is defined by Sommese in
\cite[Definition~1.3]{So}: $E$ is called \emph{$1$-ample}, if $\sO_{\P(E^*)}(1)$ is semi-ample and suppose  for some $k>0$, $\sO_{\P(E^*)}(k)$ is globally generated, then the maximum
dimension of the fiber of the evaluation map
$$\P(E^*)\>\P\left(H^0(\P(E^*),\sO_{\P(E^*)}(k))\right)$$
is $\leq 1$.  It is also pointed out in \cite[p.~127]{SS} that, $1$-ampleness is irrelevant to the metric positivity of $E$ (cf. Theorem \ref{main0}).\\

In this paper, we investigate big vector bundles and \textbf{complex manifolds} with semi-positive tangent bundles, i.e. the tangent bundles are semi-positive in the sense of Griffiths, or
equivalently, there exist \textbf{Hermitian metrics} (not necessarily K\"ahler) with semi-positive holomorphic bisectional curvature.

 The first main result of our paper can be viewed as an ample
 analogue of Kawamata-Reid-Shokurov base point free theorem for
 tangent bundles.

\btheorem\label{main1} Let $(X,\omega)$ be a compact K\"ahler manifold with semi-positive holomorphic bisectional curvature. Then the following statements are equivalent \bd \item $X$ is
Fano;

\item  the tangent bundle $TX$ is big;

\item the anti-canonical line bundle $K_X^{-1}=\det(TX)$ is big;

\item $c_1^n(X)>0$.
\ed \etheorem

\noindent  One may wonder whether  similar results hold for abstract
vector bundles. Unfortunately, there exists a vector bundle  $E$
which is semi-positive in the sense of Griffiths, and $\det(E)$ is
ample ( in particular, $\det(E)$ is big), but $E$ is not a big
vector bundle. Indeed, one can see clearly that the underlying
manifold structure of the tangent bundle is essentially used in the
proof of Theorem \ref{main1}.

 \bexample
Let $E=T\P^2\ts \sO_{\P^2}(-1)$ be the hyperplane bundle of $\P^2$.
It is easy to see that $E$ is semi-ample and semi-positive in the
sense of Griffiths and $\det E=\sO_{\P^2}(1)$ is ample. However, $E$
is not a big vector bundle since the second Segre number
$s_2(E)=c_1^2(E)-c_2(E)=0$ (for more details, see Example
\ref{counter}).\eexample

\noindent For abstract vector bundles, we obtain
 \bproposition\label{bundle0} Let $E$ be a nef vector
bundle over a compact K\"ahler manifold $X$. If $E$ is  a big vector
bundle, then $\det(E)$ is a big line bundle. \eproposition

\bcorollary If $X$ is a compact K\"ahler manifold with nef and big
tangent bundle, then $X$ is Fano. \ecorollary

As an application of Theorem \ref{main1}, we can \emph{classify}
compact K\"ahler manifolds with semi-positive holomorphic
bisectional curvature and big tangent bundles.

\btheorem\label{main0} Let $(X,\omega)$ be a compact K\"ahler
manifold with semi-positive holomorphic bisectional curvature.
Suppose
  $TX$ is a  big vector bundle. Then there exist non-negative numbers $k, N_1,\cdots, N_\ell$ and irreducible compact Hermitian symmetric spaces
   $M_1,\cdots, M_k$ of rank $\geq 2$ such that  $(X,\omega)$ is
 isometrically biholomorphic to
 \beq (\P^{N_1},\omega_1)\times\cdots \times (\P^{N_\ell},\omega_\ell)\times (M_1, \eta_1)\times\cdots\times (M_k,\eta_k) \eeq
where $\omega_i$, $1\leq i\leq \ell$, is a K\"ahler metric on
$\P^{N_i}$ with semi-positive holomorphic bisectional curvature and
$\eta_1,\cdots, \eta_k$ are the canonical metrics on $M_1,\cdots,
M_k$. \etheorem

 Note that, by Theorem \ref{SW}, the Fano manifold $\P(T^*\P^2)$
has nef and big tangent bundle.  On the other hand, it does not
admit any smooth \textbf{K\"ahler metric} with semi-positive
holomorphic bisectional curvature according to Theorem \ref{main0}
or Mok's uniformization theorem. However, it is still not clear
whether the tangent bundle of $\P(T^*\P^2)$ is semi-positive in the
sense of Griffiths, or equivalently, whether $\P(T^*\P^2)$ has a
smooth \textbf{Hermitian metric} with semi-positive holomorphic
bisectional curvature. According to a weaker version of a conjecture
of Griffiths (e.g. Remark \ref{Griffiths}), $\P(T^*\P^2)$ should
have a Hermitian metric with Griffiths semi-positive curvature since
the tangent bundle of $\P(T^*\P^2)$ is semi-ample. As motivated by
this question, we investigate complex manifolds with semi-positive
tangent bundles.


\btheorem\label{thm1}  Let $(X,\omega)$ be a compact Hermitian manifold with semi-positive holomorphic bisectional curvature, then \bd\item $X$ has Kodaira dimension $-\infty$; or
\item $X$ is a complex parallelizable manifold. \ed \etheorem

\noindent

\noindent We  also classify compact complex surfaces with semi-positive tangent bundles based on  results in \cite{DPS} (see also \cite{CP,Yau74}). In this classification, we only assume the
abstract vector bundle $TX$ is semi-positive in the sense of Griffiths, or equivalently, $X$ has a smooth \textbf{Hermitian metric} with semi-positive holomorphic bisectional curvature.
Hence, even if the ambient manifold is K\"ahler or projective, Mok's result can not be applied.

\newpage
\btheorem\label{thm3} Let $X$ be a compact K\"ahler surface. If $TX$ is (Hermitian) semi-positive, then $X$ is one of the following: \bd
\item $X$ is a torus;
\item $X$ is $\P^2$;
\item $X$ is $\P^1\times \P^1$;
\item $X$ is  a ruled surface over an elliptic curve, and $X$ is covered by $\C\times\P^1$.
\ed \etheorem \noindent We need to point out that it should be a coincidence that we get the same classification as in \cite{HS} where they considered K\"ahler metrics with semi-positive
holomorphic bisectional curvature. As explained in the previous paragraphs, it is still unclear whether one can derive the same classification in higher dimensional cases. In particular, we
would like to know whether one can get the same results as in Theorem \ref{main0} if the K\"ahler metric is replaced by a Hermitian metric.

 For
non-K\"ahler surfaces, we obtain

\btheorem\label{thm2} Let $(X,\omega)$ be a compact non-K\"ahler surface with semi-positive holomorphic bisectional curvature. Then $X$ is a Hopf surface.
 \etheorem

\noindent We also construct explicit Hermitian metrics with semi-positive curvature on Hopf surface $H_{a,b}$ (cf. \cite[Proposition~6.3]{DPS}). \bproposition\label{hopf} On every Hopf
surface $H_{a,b}$, there exists a Gauduchon metric with semi-positive holomorphic bisectional curvature. \eproposition

 \noindent For complex Calabi-Yau manifolds, i.e. complex manifolds with $c_1(X)=0$, we have
\bcorollary\label{CY0} Let $X$ be a complex Calabi-Yau manifold
 in the Fujiki class $\mathscr C$ (class
of manifolds bimeromorphic to K\"ahler manifolds). Suppose $X$ has a Hermitian metric $\omega$ with semi-positive holomorphic bisectional curvature, then $X$ is a torus. \ecorollary

\noindent By comparing Corollary \ref{CY0} with Proposition \ref{hopf}, we see that the Fujiki class condition in Corollary \ref{CY0} is necessary since every Hopf surface $H_{a,b}$ is a
Calabi-Yau manifold with semi-positive tangent bundle.\\

By using Theorem \ref{thm1} and the positivity of direct image
sheaves (Theorem \ref{directimage}) over complex manifolds (possibly
non-K\"ahler), we obtain new examples on K\"ahler and non-K\"ahler
manifolds whose tangent bundles  are \emph{nef but not
semi-positive}. To the best of our knowledge, it is also the first
example in the manifold setting (cf. \cite[Example~1.7]{DPS}).

\bcorollary\label{example} Let $X$ be a Kodaira surface or a hyperelliptic surface.

 \bd\item The  tangent bundle $TX$  is nef
but not semi-positive (in the sense of Griffiths);

\item The anti-canonical line
bundle of  $\P(T^*X)$ is nef,  but neither   semi-positive nor big.

 \ed\ecorollary
\noindent Hence, for any dimension $n\geq 2$, there exist K\"ahler and non-K\"ahler manifolds with nef but not semi-positive tangent bundles.\\

Finally, we investigate compact complex manifolds, of the form $\P(T^*X)$, whose tangent bundles are nef. It is well-known that $\P(T^*\P^n)$ is homogeneous, and its tangent bundle is nef.
We obtain a similar converse statement and yield another characterization of $\P^n$.

\bproposition\label{PTX1} Let $X$ be a Fano manifold of complex dimension $n$. Suppose $\P(T^*X)$ has nef tangent bundle, then $X\cong \P^n$. \eproposition

\noindent In particular, for complex $3$-folds, we have the following classification.

 \btheorem\label{PTX2} For a complex $3$-fold $\P(T^*X)$, if
 $\P(T^*X)$ has nef tangent bundle, then $X$ is
exactly one of the following: \bd \item $X\cong \P^2$;\item $X\cong \T^2$, a flat torus;

\item $X$ is a hyperelliptic surface;

\item $X$ is a Kodaira surface;

\item $X$ is a Hopf surface.
\ed

\etheorem

The paper is organized as follows: In Section \ref{basic}, we introduce several basic terminologies which will be used frequently in the paper. In Section \ref{direct}, we study the
positivity of direct image sheaves over complex manifolds (possibly non-K\"ahler). In Section \ref{kahlerbig}, we investigate compact K\"ahler manifolds with big tangent bundles and prove
 Theorem \ref{main1}, Proposition \ref{bundle0} and Theorem \ref{main0},. In Section \ref{semi}, we study compact complex manifolds with semi-positive tangent bundles and establish Theorem
\ref{thm1}, Theorem \ref{thm3}, Theorem \ref{thm2}, Proposition
\ref{hopf}, Corollary \ref{CY0} and Corollary \ref{example}. In
Section \ref{PTX}, we discuss complex manifolds of the form
$\P(T^*X)$ and prove Proposition \ref{PTX1} and Theorem \ref{PTX2}.
In the Appendix \ref{appendix}, we include some straightforward
computations on Hopf manifolds for the reader's convenience.

\bremark For compact K\"ahler manifolds with semi-negative holomorphic bisectional curvature, there are similar uniformization theorems as
 Mok's result. We refer to \cite{WZ}, \cite{L} and the references
 therein. We have obtained a number of  results for compact complex manifolds with semi-negative tangent bundles, which will
 appear in \cite{Yang}.

\eremark

\noindent\textbf{Acknowledgement.} The author would like to thank Professor
 K.-F. Liu, L.-H. Shen,  V. Tosatti, B. Weinkove, S.-T. Yau, and Y. Yuan  for many
valuable discussions. The author would also like to thank Professor T. Peternell for answering his question, which leads to the current version of Proposition \ref{surface}.

\newpage

\section{Background materials}\label{basic}

Let $E$ be a holomorphic vector bundle over a compact complex manifold $X$ and $h$ a Hermitian metric on $E$. There exists a unique connection $\nabla$ which is compatible with the
 metric $h$ and the complex structure on $E$. It is called the Chern connection of $(E,h)$. Let $\{z^i\}_{i=1}^n$ be  local holomorphic coordinates
  on $X$ and  $\{e_\alpha\}_{\alpha=1}^r$ be a local frame
 of $E$. The curvature tensor $R^\nabla\in \Gamma(X,\Lambda^2T^*X\ts E^*\ts E)$ has components \beq R_{i\bar j\alpha\bar\beta}= -\frac{\p^2
h_{\alpha\bar \beta}}{\p z^i\p\bar z^j}+h^{\gamma\bar \delta}\frac{\p h_{\alpha \bar \delta}}{\p z^i}\frac{\p h_{\gamma\bar\beta}}{\p \bar z^j}\eeq Here and henceforth we sometimes adopt the
Einstein convention for summation.

 \bdefinition
A Hermitian holomorphic vector bundle $(E,h)$ is called positive (resp. semi-positive) in the sense of Griffiths if $$ R_{i\bar j \alpha\bar\beta} u^i\bar u^j v^\alpha\bar v^\beta>0 \qtq{(
resp. $\geq 0$ )}$$ for nonzero vectors $u=(u^1,\cdots, u^n)$ and $v=(v^1,\cdots, v^r)$ where $n=\dim_\C X$ and $r$ is the rank of $E$. $(E,h)$ is called Nakano positive (resp. Nakano
semi-positive) if
$$ R_{i\bar j
\alpha\bar\beta} u^{i\alpha}\bar u^{j\beta}>0\qtq{( resp. $\geq 0$)}
$$
for nonzero vector $u=(u^{i\alpha})\in \C^{nr}$.

 \edefinition

 In
particular, if $(X,\omega_g)$ is a  Hermitian manifold, $(T^{1,0}M,\omega_g)$ has  Chern curvature components \beq R_{i\bar j k\bar \ell}=-\frac{\p^2g_{k\bar \ell}}{\p z^i\p\bar z^j}+g^{p\bar
q}\frac{\p g_{k\bar q}}{\p z^i}\frac{\p g_{p\bar \ell}}{\p\bar z^j}.\eeq  The (first) Chern-Ricci form $Ric(\omega_g)$ of $(X,\omega_g)$ has components
$$R_{i\bar j}=g^{k\bar \ell}R_{i\bar jk\bar \ell}=-\frac{\p^2\log\det(g)}{\p z^i\p\bar z^j}$$
and it is well-known that the Chern-Ricci form represents the first Chern class of the complex manifold $X$ (up to a factor $2\pi$).

\bdefinition Let $(X,\omega)$ be a compact Hermitian manifold. $(X,\omega)$ has positive (resp. semi-positive)  holomorphic \emph{ bisectional} curvature, if for any nonzero vector
$\xi=(\xi^1,\cdots, \xi^n)$ and $\eta=(\eta^1,\cdots, \eta^n)$,
$$R_{i\bar j k\bar \ell}\xi^i\bar\xi^j\eta^k\bar\eta^\ell>0\ \ \ \text{(resp. $\geq 0$)}.$$
$(X,\omega)$ has positive (resp. semi-positive) holomorphic
 \emph{sectional} curvature, if for any nonzero vector $\xi=(\xi^1,\cdots,
\xi^n)$
$$R_{i\bar j k\bar \ell}\xi^i\bar\xi^j\xi^k\bar\xi^\ell>0\ \ \ \text{(resp. $\geq 0$)}.$$
 \edefinition

\bdefinition Let $(X,\omega)$ be a Hermitian manifold, $L\>X$ a holomorphic line bundle and $E\>X$ a holomorphic vector bundle. Let $\sO_{\P(E^*)}(1)$ be the tautological line bundle of the
projective bundle $\P(E^*)$ over $X$. \bd\item $L$ is said to be \emph{positive} (resp. \emph{semi-positive}) if there exists a smooth Hermitian metric $h$ on $L$ such that the curvature form
$R=-\sq\p\bp\log h$ is a positive (resp. semi-positive) $(1,1)$-form. The vector bundle $E$ is called \emph{ample} (resp. \emph{semi-ample}) if $\sO_{\P(E^*)}(1)$  is a positive (resp.
semi-positive) line bundle.

\item $L$ is said to be \emph{nef} ( or numerically effective), if for any  $\eps>0$, there exists a
smooth Hermitian metric $h$ on $L$ such that the curvature of $(L,h)$ satisfies $ -\sq\p\bp\log h\geq -\eps \omega.$ The vector bundle $E$ is called \emph{nef} if $\sO_{\P(E^*)}(1)$  is a nef
line bundle.

\item $L$ is said to be \emph{big}, if there exists a
(possibly) singular Hermitian metric $h$ on $L$ such that the curvature  of $(L,h)$ satisfies $R=-\sq\p\bp\log h\geq \eps \omega$
 in the sense of current for some $\eps> 0$. The vector bundle $E$ is called \emph{big}, if  $\sO_{\P(E^*)}(1)$  is big.

 \ed\edefinition

\bdefinition Let $X$ be a compact complex manifold and $L\>X$ be a line bundle. The Kodaira dimension $\kappa(L)$ of $L$ is defined to be
$$\kappa(L):=\limsup_{m\>+\infty} \frac{\log \dim_\C
H^0(X,L^{\ts m})}{\log m}$$ and the \emph{Kodaira dimension} $\kappa(X)$ of $X$ is defined as $ \kappa(X):=\kappa(K_X)$ where the logarithm of zero is defined to be $-\infty$. \edefinition

\noindent By Riemann-Roch, it is easy to see that $E$ is a big vector bundle if and only if there are $c_0>0$ and $k_0\geq 0$ such that \beq h^0(X,S^kE)\geq c_0 k^{n+r-1}\eeq for all $k\geq
k_0$ where $\dim_\C X=n$ and $rk(E)=r$. Indeed, let $Y=\P(E^*)$ and $\sO_Y(1)$ be the tautological line bundle of $Y$, then we have \beq h^0(X,S^kE)=h^0(Y,\sO_Y(k))\geq c_0k^{n+r-1}\eeq where
$\dim_\C Y=n+r-1$. Hence, $E$ is big if and only if $\sO_{\P(E^*)}(1)$ is big, if and only if
$$\kappa(\sO_{\P(E^*)}(1))=\dim_\C \P(E^*).$$

\noindent The following well-known result will be used frequently in
the paper. \blemma Let $L$ be a line bundle over a compact K\"ahler
manifold $X$. Suppose $L$ is nef, then $L$ is big if and only if the
top self-intersection number $c^n_1(L)>0$ where $n=\dim X$. \elemma

\section{Positivity of direct image sheaves over  complex
 manifolds}\label{direct}

Let $\sX$ be a compact complex manifold of complex dimension $m+n$, and $S$ a smooth complex manifold (\emph{possibly non-K\"ahler}) with dimension $m$. Let $\pi:\sX\>S$ be a smooth proper
submersion such that for any $s\in S$, $X_s:=\pi^{-1}(\{s\})$ is a compact K\"ahler manifold with dimension $n$. Suppose for any $s\in S$, there exists an open neighborhood $U_s$ of $s$ and a
smooth $(1,1)$ form $\omega$ on $\pi^{-1}(U_s)$ such that $\omega_p=\omega|_{X_p}$ is a smooth K\"ahler form on $X_p$ for any $p\in U_s$.
 Let $(\sE,h^{\sE})\>\sX$ be a
Hermitian holomorphic vector bundle. In the following, we adopt the setting in \cite[Section ~4]{Bo} (see also \cite[Section~2.3]{LYJDG}). Consider the space of holomorphic $\sE$-valued
$(n,0)$-forms on $X_s$,
 $$E_s:=H^0(X_s,\sE_s\ts K_{X_s})\cong H^{n,0}(X_s,\sE_s)$$
where $\sE_s=\sE|_{X_s}$. Here, we assume all $E_s$ have the same dimension. With a natural holomorphic structure,
$$E=\bigcup_{s\in S}\{s\}\times E_s$$ is isomorphic to the  direct
image sheaf $\pi_*(K_{\sX/S}\ts \sE)$ if $\sE$ has certain positive property. \btheorem\label{directimage}If $(\sE,h^{\sE})$ is positive (resp. semi-positive) in the sense of Nakano, then
$\pi_*(K_{\sX/S}\ts \sE)$ is positive (resp. semi-positive) in the sense of Nakano. \etheorem \bremark When the total space $\sX$ is K\"ahler and $\sE$ is a line bundle, Theorem
\ref{directimage} is originally proved by Berndtsson in \cite[Theorem~1.2]{Bo}. When $(\sE,h^{\sE})$ is a Nakano semi-positive vector bundle, Theorem \ref{directimage} is a special case of
\cite[Theorem~1.1]{MT}. \eremark

It is not hard to see that the positivity of the direct image sheaves does not  depend on the base manifold $S$. It still works for non-K\"ahler $S$. We give a sketched proof of Theorem
\ref{directimage} for reader's convenience.
 Let $h^{\sE}$ be a smooth Nakano semi-positive metrics on
$\sE$. For any local smooth section $u$ of $\pi_*(K_{\sX/S}\ts \sE)$, there is a representative $\textbf{u}$ of $u$, a local holomorphic $\sE$-valued $(n,0)$ form on $\sX$, then we define the
Hodge metric on $\pi_*(K_{\sX/S}\ts \sE)$ by using the sesquilinear pairing \beq |u|^2=\sq \int_{X_s}\{\textbf{u},\textbf{u}\}.\eeq Note that we do not specify any metric on $\sX$ or $S$.
Since $\sX\>S$ has K\"ahler fibers, we can use similar methods as in \cite{Bo,LYJDG} to compute the curvature of the Hodge metric. To obtain the positivity of the Hodge metric, the key
ingredient is to find primitive representatives on the K\"ahler fiber $X_s$ (e.g. \cite[Lemma~4.3]{Bo} or \cite[Theorem~3.10]{LYJDG}). Since all computations are local, i.e. on an open subset
$\pi^{-1}(U)$ of $\sX$ where $U$ is an open subset of $S$, the computations do not depend on the property of base manifold $S$. In particular, all computations in \cite{LYJDG} and all results
(e.g. \cite[Theorem~1.1 and Theorem~1.6]{LYJDG}) still work for non-K\"ahler base manifold $S$. Note that, if $(\sE,h^{\sE})$ is only semi-positive, $S$ can be a non-K\"ahler manifold.

\bcorollary\label{directimage2} Let $X$ be a compact complex manifold (possibly non-K\"ahler) and $E\>X$ be a holomorphic vector bundle of rank $r$.

\bd \item If $\sO_{\P(E^*)}(1)$ is semi-positive, then $\emph{S}^kE\ts \det(E)$ is Nakano semi-positive.

 \item If $\det E$ is a holomorphic torsion, i.e.  $(\det E)^{ k}=\sO_X$ for some $k\in \N^+$, then $E$ is Nakano semi-positive if and only
if $\sO_{\P(E^*)}(1)$ is semi-positive. \ed \ecorollary

\bproof Let $Y=\P(E^*)$, $L=\sO_{\P(E^*)}(1)$ and $\pi:\P(E^*)\>X$ be the canonical projection.  $(1).$  By the adjunction formula \cite[p.~89]{Laza1}, we have \beq K_Y=L^{-r}\ts \pi^*(K_X\ts
\det (E)),\label{adjunction}\eeq and \beq K_{Y/X}=L^{-r}\ts \pi^*(\det(E)).\label{adjunction2}\eeq Therefore, $$\pi_*(K_{Y/X}\ts L^{r+k})=\pi_*(L^k\ts \pi^*(\det E))=\emph{S}^kE\ts \det E.$$
By Theorem \ref{directimage}, we deduce $S^kE\ts \det E$ is semi-positive in the sense of Nakano if $L$ is semi-positive.\\ $(2).$ Suppose $\det E$ is a holomorphic torsion with $(\det
E)^m=\sO_X$, then there exists a flat Hermitian metric on $\det E$ and also on $\det E^*$. If $L$ is semi-positive, $$\tilde L=L^{r+1}\ts \pi^*(\det E^*)$$ is semi-positive. By formula
(\ref{adjunction2}) and Theorem \ref{directimage}, we know
$$\pi_*(K_{Y/X}\ts \tilde L)=\pi_*(L)=E$$ is semi-positive in the
sense of Nakano.

 On the other hand, if
$(E,h)$ is semi-positive, then the induced Hermitian metric on $L$ has semi-positive curvature (e.g. formula (\ref{curvature})).\eproof

\bremark\label{Griffiths} Griffiths conjectured in \cite{G} that $E$ is Griffiths positive  if (and only if) the tautological line bundle $\sO_{\P(E^*)}(1)$ is positive. It is also not known
in the semi-positive setting, i.e. whether there exists a Griffiths semi-positive metric on $E$ when $\sO_{\P(E^*)}(1)$ is semi-positive. \eremark

\section{K\"ahler manifolds with big tangent
bundles}\label{kahlerbig}

In this section, we prove  Theorem \ref{main0} and Theorem \ref{main1}. We begin with an algebraic curvature relation on a K\"ahler manifold $(X,\omega)$. At a given point $p\in X$, the
minimum holomorphic sectional curvature is defined to be
$$\min_{W\in T^{1,0}_pX,|W|=1}H(W),$$ where $H(W):=R(W,\bar W,
W,\bar W)$. Since $X$ is of finite dimension, the minimum can be attained. \blemma\label{linear} Let $(X,\omega)$ be a compact K\"ahler manifold. Let $e_1\in T_p^{1,0}X$ be a unit vector
which minimizes the holomorphic sectional curvature of $\omega$ at point $p$, then \beq 2R(e_1,\bar e_1,W,\bar W)\geq \left(1+|\la W,e_1\ra|^2\right) R(e_1,\bar e_1,e_1,\bar e_1) \eeq for
every unit vector $W\in T^{1,0}_pX$.

 \elemma \bproof
 Let $e_2\in T^{1,0}_pX$ be any unit vector  orthogonal
to $e_1$. Let $$f_1(\theta)= H(\cos(\theta )e_1+\sin(\theta) e_2),\ \ \ \ \theta\in \R.$$ Then we have \be f_1(\theta)&=&R(\cos(\theta )e_1+\sin(\theta) e_2, \bar{\cos(\theta
)e_1+\sin(\theta) e_2},\cos(\theta )e_1+\sin(\theta) e_2,\bar{\cos(\theta )e_1+\sin(\theta) e_2})\\&=&\cos^4(\theta)R_{1\bar 1 1\bar 1}+\sin^4(\theta)R_{2\bar 2 2 \bar
2}\\&&+2\sin(\theta)\cos^3(\theta)\left[R_{1\bar 1 1\bar 2}+R_{2\bar 1 1\bar 1}\right]+2\cos(\theta)\sin^3(\theta)\left[R_{1\bar 2 2\bar 2}+R_{2\bar 1 2\bar 2}\right]\\
&&+\sin^2(\theta)\cos^2(\theta)[4R_{1\bar 1 2\bar 2}+R_{1\bar 21\bar 2}+R_{2\bar 1 2\bar 1}].\ee Since $f_1(\theta)\geq R_{1\bar 1 1\bar 1}$ for all $\theta\in \R$ and  $f_1(0)=R_{1\bar
11\bar 1}$, we have
$$f_1'(0)=0\qtq{ and} f_1''(0)\geq 0.$$ By a straightforward
computation, we obtain \beq f_1'(0)=2(R_{1\bar 1 1\bar 2}+R_{2\bar 1 1\bar 1})=0,\ \ \ \ f_1''(0)=2\left(4R_{1\bar 1 2\bar 2}+R_{1\bar 21\bar 2}+R_{2\bar 1 2\bar 1}\right)- 4R_{1\bar 1 1\bar
1}\geq 0. \label{a}\eeq Similarly, if we set $f_2(\theta)=H(\cos(\theta )e_1+\sq \sin(\theta) e_2)$, then \be f_2(\theta)&=&\cos^4(\theta)R_{1\bar 1 1\bar 1}+\sin^4(\theta)R_{2\bar 2 2 \bar
2}\\&&+2\sq \sin(\theta)\cos^3(\theta)\left[- R_{1\bar 1 1\bar 2}+ R_{2\bar 1 1\bar 1}\right]+2\sq\cos(\theta)\sin^3(\theta)\left[- R_{1\bar 2 2\bar 2}+ R_{2\bar 1 2\bar 2}\right]\\
&&+\sin^2(\theta)\cos^2(\theta)[4R_{1\bar 1 2\bar 2}-R_{1\bar 21\bar 2}-R_{2\bar 1 2\bar 1}].\ee

\noindent  From  $f'_2(0)=0$ and $f''_2(0)\geq 0$, one can see \beq -R_{1\bar 1 1\bar 2}+R_{2\bar 1 1\bar 1}=0,\ \ \ \ \ \ \ 2\left(4R_{1\bar 1 2\bar 2}-R_{1\bar 21\bar 2}-R_{2\bar 1 2\bar
1}\right)-4R_{1\bar 1 1\bar 1}\geq 0.\label{b}\eeq Hence,  from (\ref{a}) and (\ref{b}), we obtain \beq R_{1\bar 1 1\bar 2}=R_{1\bar 1 2\bar 1}=0, \qtq{and} 2R_{1\bar 1 2\bar 2}\geq R_{1\bar
1 1 \bar 1}.\label{bbb}\eeq For an arbitrary unit vector $W\in T^{1,0}_pX$, if $W$ is parallel to $e_1$, i.e. $W=\lambda e_1$ with $|\lambda|=1$,
$$2R(e_1,\bar e_1, W,\bar W)=2R(e_1,\bar e_1, e_1, \bar e_1).$$
Suppose $W$ is not parallel to $e_1$. Let $e_2$ be the unit vector
$$e_2=\frac{W-\la W,e_1\ra e_1}{|W-\la W,e_1\ra e_1|}.$$ Then $e_2$
is a unit vector orthogonal to $e_1$ and
$$W=ae_1+be_2,\ \ \ \ a=\la W,e_1\ra,\ \ \ b=|W-\la W,e_1\ra|,\ \ \
|a|^2+|b|^2=1.$$ Hence $$2R(e_1,\bar{e}_1,W,\bar W)=2|a|^2R_{1\bar 1 1\bar 1}+2|b|^2 R_{1\bar 12\bar 2},$$ since we have $R_{1\bar 1 1\bar 2}=R_{1\bar 1 2\bar 1}=0$ by (\ref{bbb}). By
(\ref{bbb}) again,
$$2R(e_1,\bar{e}_1,W,\bar W)\geq (2|a|^2+|b|^2)R_{1\bar 1 1\bar 1}=(1+|a|^2) R_{1\bar 1 1\bar 1}$$
which completes the proof of the lemma.
 \eproof

\noindent By using similar methods, one has

\blemma\label{linear1}  Let $e_n\in T_p^{1,0}X$ be a unit vector which maximizes the holomorphic sectional curvature  at point $p$, then \beq 2R(e_n,\bar e_n,W,\bar W)\leq \left(1+|\la W,
e_n\ra|^2\right) R(e_n,\bar e_n,e_n,\bar e_n) \eeq for every unit vector $W\in T^{1,0}_pX$. \elemma

\bremark A special case of Lemma \ref{linear1}--when $W$ is orthogonal to $e_n$--is well-known (e.g. \cite[p.~312]{Goldberg}, \cite[p.~136]{Brendle}).  When the holomorphic sectional
curvature is strictly negative at point $p$, one has $2R(e_n,\bar e_n,W,\bar W)\leq R(e_n,\bar e_n,e_n,\bar e_n)$,  which is firstly obtained in \cite[Lemma~1.4]{BKT}. In the proofs of Lemma
\ref{linear} and Lemma \ref{linear1}, we refine the methods in \cite{Goldberg} and \cite{Brendle}.

\eremark

\btheorem\label{thm00} Let $(X,\omega)$ be a compact K\"ahler manifold with semi-positive holomorphic bisectional curvature. Then the following statements are equivalent \bd \item $X$ is
Fano;

\item  the tangent bundle $TX$ is big;

\item the anti-canonical line bundle $K_X^{-1}=\det(TX)$ is big;

\item $c_1^n(X)>0$.
\ed \etheorem

\bproof $(1)\Longrightarrow(2)$. Let $E=TX$ and $L=\sO_{\P(E^*)}(1)$ the tautological line bundle over the projective bundle $\P(E^*)$. Let's recall the general setting when $(E,h^E)$ is an
arbitrary Hermitian holomorphic vector bundle (e.g.\cite{G,Demailly,LYJAG}). Let $(e_1,\cdots, e_n)$ be the local holomorphic frame with respect to a given trivialization on $E$ and the dual
frame on $E^*$ is denoted by $(e^1,\cdots, e^n)$. The corresponding holomorphic coordinates on $E^*$ are denoted by $(W_1,\cdots, W_n)$.  There is a local section $e_{L^*}$ of $L^*$ defined
by $$ e_{L^*}=\sum_{\alpha=1}^n W_\alpha e^\alpha.$$ Its dual section is denoted by $e_L$. Let $h^E$ be a fixed Hermitian metric on $E$ and $h^L$ the induced quotient metric by the morphism
$(\pi^*E,\pi^*h^E)\>L$. %

If $\left(h_{\alpha\bar\beta}\right)$ is the matrix representation of $h^E$ with respect to the basis $\{e_\alpha\}_{\alpha=1}^n$, then $h^L$ can be written as \beq
h^L=\frac{1}{h^{L^*}(e_{L^*},e_{L^*})}=\frac{1}{\sum h^{\alpha\bar\beta}W_\alpha\bar W_\beta}.\label{inducedmetric} \eeq
 The curvature of $(L,h^L)$ is \beq
R^{h^L}=-\sq\p\bp\log h^L=\sq\p\bp\log\left(\sum h^{\alpha\bar\beta}W_\alpha\bar W_\beta\right) \label{inducedcurvature}\eeq where $\p$ and $\bp$ are operators on the total space $\P(E^*)$.
 We fix a point $Q\in \P(E^*)$, then there exist local
holomorphic coordinates
 $(z^1,\cdots, z^n)$ centered at point $p=\pi(Q)$ and local holomorphic basis $\{e_1,\cdots, e_n\}$ of $E$ around $p\in X$ such that
 \beq h_{\alpha\bar\beta}=\delta_{\alpha\bar\beta}-R_{i\bar j \alpha\bar\beta}z^i\bar z^j+O(|z|^3). \label{normal}\eeq
Without loss of generality, we assume $Q$  is the point $(0,\cdots,
0,[a_1,\cdots, a_n])$ with $a_n=1$. On the chart $U=\{W_n=1\}$ of
the fiber $\P^{n-1}$, we set $w^A=W_A$ for $A=1,\cdots, n-1$. By
formula (\ref{inducedcurvature}) and (\ref{normal}), we obtain the
well-known formula (e.g.\cite[Proposition~2.5]{LYJAG}) \beq
R^{h^L}(Q)=\sq\left(\sum_{\alpha,\beta=1}^n R_{i\bar
j\alpha\bar\beta}\frac{a_\beta \bar a_\alpha}{|a|^2}dz^i\wedge d\bar
z^j+\sum_{A,B=1}^{n-1}\left(1-\frac{a_B\bar
a_A}{|a|^2}\right)dw^A\wedge d\bar w^B\right)\label{curvature} \eeq
where $|a|^2=\sum\limits_{\alpha=1}^n|a_\alpha|^2$.

 Since $(X,\omega)$ is a K\"ahler manifold with semi-positive
 holomorphic bisectional curvature, the Ricci curvature $Ric(\omega)$ of
 $\omega$ is also semi-positive. On the other hand, since $X$ is Fano,
 we have
 $$\int_X \left(Ric(\omega)\right)^n>0.$$
Therefore, $Ric(\omega)$ must be strictly positive at some point $p\in X$. Then by a result of Mok \cite[Proposition~1.1]{Mok}, there exists a K\"ahler metric $\hat \omega$ such that $\hat
\omega$ has semi-positive holomorphic bisectional curvature, strictly positive holomorphic sectional curvature and strictly positive Ricci curvature. Indeed, let \beq
\begin{cases} \frac{\p \omega_t }{\p t}=-Ric(\omega_t),\\
\omega_0=\omega
\end{cases}
\eeq be the K\"ahler-Ricci flow with initial metric $\omega$, then
we can take $\omega_t$ as $\hat \omega$ for some small positive $t$
satisfying $[\omega]-tc_1(X)>0$. Let $\hat R$ be the corresponding
curvature operator of $\hat \omega$.
 We choose normal coordinates $\{z^1,\cdots, z^n\}$ centered at point $p$ such that
 $\{e_i=\frac{\p}{\p z^i}\}_{i=1}^n$ is  the normal frame of $(E,\hat \omega)=(TX,\hat\omega)$. Let $K\in T^{1,0}_pX$ be
 a unit vector which minimizes the holomorphic sectional curvature
 of $\hat \omega$
 at point $p\in X$. In particular, we have $\hat R(K,\bar K,K,\bar K)>0$. Hence there exists a unit vector $a=(a_1,\cdots, a_n)\in \C^n$  such that
 \beq K=a_1e_1+\cdots +a_n e_n.\eeq
Without loss of generality, we assume $a_n\neq 0$. By Lemma \ref{linear}, for any unit vector $W=\sum_{} b_i e_i \in T_p^{1,0}X$, we have \beq \hat R(K,\bar K, W,\bar W)\geq \frac{1}{2} \hat
R(K,\bar K,K,\bar K)>0.\eeq That is \beq\sum_{i,j,k,\ell} \hat R_{i\bar j k\bar \ell}a_k\bar a_\ell b_i\bar b_j >0\label{key}\eeq for every unit vector $b=(b_1,\cdots, b_n)$ in $\C^n$.
 Then at
point $Q\in Y=\P(T^*X)$ with coordinates
 $$(0,\cdots, 0, [\bar a_1,\cdots, \bar a_n])=\left(0,\cdots, 0,\left[\frac{\bar a_1}{\bar a_n},\cdots, \frac{\bar a_{n-1}}{\bar a_n}, 1\right]\right),$$ we obtain
 \beq R^{h^L}(Q)=\sq\left(\sum_{k,\ell=1}^n \hat R_{i\bar
jk\bar\ell}a_k \bar a_\ell dz^i\wedge d\bar
z^j+\sum_{A,B=1}^{n-1}\left(1- a_A\bar a_B\right)dw^A\wedge d\bar
w^B\right)\label{key2}\eeq is a strictly positive $(1,1)$ form at
point $Q\in Y$ according to (\ref{key}). Here, we also use equation
(\ref{curvature}) and the fact that $\sum_{i=1}^n|a_i|^2=1$. By
continuity, $(L,h^L)$ is positive at a small neighborhood of $Q$.
Since we already know $c_1(L)\geq 0$, and so
$$\int_Y c^{2n-1}_1(L)>0.$$ Hence $L$ is a big line bundle by
Siu-Demailly's solution to the Grauert-Riemenschneider conjecture
 (\cite{Siu,Dem85}).
 In particular, the  tangent bundle $TX$ is big.
\\
$(3)\Longleftrightarrow(4).$ Since $K^{-1}_X$ is semi-positive and
in particular it is nef,
it is well-known  that they are equivalent.\\
$(4)\Longrightarrow(1).$ This part is well-known (e.g.
\cite[Theorem~4.2]{DPS}), we include a sketch for reader's
convenience. Since $TX$ is nef, and so is $K_X^{-1}=\det(TX)$. If
$c_1^n(X)>0$, we know $K_X^{-1}$ is nef and big. Hence $X$ is
K\"ahler and Moishezon, and so it is projective. By
Kawamata-Reid-Shokurov base point free theorem (e.g.
\cite[Theorem~3.3]{KM}), $K_X^{-1}$ is semi-ample, i.e. $K^{-m}_X$
is generated by global sections for some large $m$. Let
$\phi:X\dashrightarrow Y$ be the birational map defined by
$|K_X^{-m}|$. If $K^{-1}_X$ is not ample, then there exists a
rational curve $C$ contracted by $\phi$. Since $TX$ is nef, $C$
deforms to cover $X$ which is a contradiction.\\
$(2)\Longrightarrow(4).$ Since $TX$ is nef, $K_X^{-1}$ is also nef.
In particular, we have $c_1^n(X)=c^n_1(TX)\geq 0$. If $c_1^n(X)=0$,
then all Chern numbers of $X$ are zero (\cite[Corollary~2.7]{DPS}).
On the other hand,  since the signed Segre number $(-1)^ns_n(TX)$ is
a combination of Chern numbers [e.g. formula (\ref{recursion})], we
deduce that
$$(-1)^ns_n(TX)=0.$$ Hence $TX$ is not big by Lemma
\ref{bigbundle}.\eproof

\btheorem\label{main00} Let $(X,\omega)$ be a compact K\"ahler manifold with semi-positive holomorphic bisectional curvature. Suppose
  $TX$ is a  big vector bundle. Then there exist non-negative numbers $k, N_1,\cdots, N_\ell$ and irreducible compact Hermitian symmetric spaces
   $M_1,\cdots, M_k$ of rank $\geq 2$ such that  $(X,\omega)$ is
 isometrically biholomorphic to
 \beq (\P^{N_1},\omega_1)\times\cdots \times (\P^{N_\ell},\omega_\ell)\times (M_1, \eta_1)\times\cdots\times (M_k,\eta_k) \eeq
where $\omega_i$, $1\leq i\leq \ell$, is a K\"ahler metric on $\P^{N_i}$ with semi-positive holomorphic bisectional curvature and $\eta_1,\cdots, \eta_k$ are the canonical metrics on
$M_1,\cdots, M_k$. \etheorem

\bproof By Theorem \ref{thm00}, $X$ is Fano. By Yau's theorem \cite{Yau}, there exists a K\"ahler metric with strictly positive Ricci curvature. Hence $\pi_1(X)$ is finite by Myers' theorem.
By Kodaira vanishing theorem, for any $q\geq 1$, $H^{0,q}(X)=H^{n,q}(X,K_X^{-1})=0$ since $K_X^{-1}$ is ample. Therefore the Euler-Poincar\'e characteristic $\chi(\sO_X)=\sum
(-1)^qh^{0,q}(X)=1$. Let $\tilde X$ be the universal cover of $X$. Suppose it is a $p$-sheet cover over $X$, where $p=|\pi_1(X)|$. So $\tilde X$ is still a Fano manifold and hence
$\chi(\sO_{\tilde X})=p\cdot \chi(\sO_X)=1$. We obtain $p=1$, i.e. $X$ is simply connected and $\tilde X=X$. By Mok's uniformization theorem (\cite{Mok}) for compact K\"ahler manifolds with
\emph{semi-positive} holomorphic bisectional curvature, $\tilde X=X$ is isometrically biholomorphic to \beq (\P^{N_1},\omega_1)\times\cdots \times (\P^{N_\ell},\omega_\ell)\times (M_1,
\eta_1)\times\cdots\times (M_k,\eta_k) \eeq where $\omega_i$, $1\leq i\leq \ell$, is a K\"ahler metric on $\P^{N_i}$ with semi-positive holomorphic bisectional curvature and $\eta_1,\cdots,
\eta_k$ are the canonical metrics on the irreducible compact Hermitian symmetric spaces $M_1,\cdots, M_k$. Note also that, all irreducible compact Hermitian symmetric spaces (with rank $\geq
2$) are Fano. \eproof

As an application of Theorem \ref{main00}, we have

 \bcorollary\label{pn} Let
$X=\P^m\times \P^n$ and $Y=\P(T^*X)$. Then \bd
\item the tangent bundle $TX$ of $X$ is nef and big;
\item the anti-canonical line bundle $K^{-1}_Y$ of $Y$ is nef,
big, semi-ample, quasi-positive but not ample;
\item
the holomorphic tangent bundle $TY$ is not nef. \ed \ecorollary \bproof $(1)$ is from Theorem \ref{main00}. $(2)$. By adjunction formula (\ref{adjunction}), $$K_Y^{-1}=\sO_{Y}(m+n)$$ where
$\sO_Y(1)$ is the tautological line bundle of the projective bundle $\P(T^*X)$. Hence, $K_Y^{-1}$ is nef and big, and so is semi-ample by Kawamata-Reid-Shokurov base point free theorem. Let
$\omega$ be the K\"ahler metric on $X=\P^m\times \P^n$ induced by the Fubini-Study metrics. It is easy to see that $\omega$ has semi-positive holomorphic bisectional curvature and strictly
positive holomorphic sectional curvature. By Lemma \ref{linear} and formula (\ref{key2}), the induced Hermitian metric on $L=\sO_Y(1)$ is quasi-positive, i.e. $\sO_{Y}(1)$ is semi-positive
and strictly positive at some point. In particular,
 $K_Y^{-1}$ is quasi-positive. However, $K_Y^{-1}$ is not ample, otherwise $\sO_Y(1)$ is ample and so it $TX$. $(3)$. If $TY$ is nef, then the nef and big line bundle $K_Y^{-1}$ is ample.
 \eproof

As motivated by Theorem \ref{thm00}, we investigate properties for abstract nef and big vector bundles. Let $c(E)$ be the total Chern class of a vector bundle $E$, i.e.
$c(E)=1+c_1(E)+\cdots+c_n(E)$. The total Segre class $s(E)$ is defined to be the inverse of the total Chern class, i.e.
$$c(E)\cdot s(E)=1$$ where $s(E)=1+s_1(E)+\cdots+s_n(E)$ and
$s_k(E)\in H^{2k}(X)$, $1\leq k\leq n$. We have the recursion formula \beq s_k(E)+s_{k-1}(E)\cdot c_1(E)+\cdots+s_1(E)\cdot c_{k-1}(E)+c_k(E)=0\label{recursion}\eeq for every $k\geq 1$. In
particular, one has \beq s_1(E)=-c_1(E),\ s_2(E)=c_1^2(E)-c_2(E),\ s_3(E)=2c_1(E)c_2(E)-c_1^3(E)-c_3(E).\eeq In particular, the top Segre class $s_n(E)$ is a polynomial of Chern classes of
degree $2n$. (Note that there is alternated sign's difference from the notations in \cite[p.~245]{G}).
 The following result is essentially well-known.

\blemma\label{bigbundle} Let $(X,\omega)$ be a compact K\"ahler manifold with complex dimension $n$. Suppose $E$ is nef vector bundle  with rank $r$, then $E$ is  big if and only if the
signed Segre number $(-1)^ns_n(E)>0$. \elemma \bproof Let $L=\sO_{\P(E^*)}(1)$ and $\pi:\P(E^*)\>X$ be the canonical projection. Since $L$ is nef,  $L$ is a big line bundle if and only if the
top self intersection number
$$c^{n+r-1}_1(L)>0.$$ On the other hand, by \cite[Proposition~5.22]{G} we have  \beq \pi_*(c_1^{n+r-1}(L))=(-1)^n
s_n(E),\eeq where $\pi_*:H^{2(n+r-1)}(\P(E^*))\>H^{2n}(X)$ is the pushforward homomorphism induced by $\pi$.
 Hence $L$ is big if and only if the signed Segre number
$(-1)^ns_n(E)$ is positive. \eproof

\bproposition\label{bundle} Let $E$ be a nef vector bundle over a compact K\"ahler manifold $X$. If $E$ is  a big vector bundle, then $\det(E)$ is a big line bundle. \eproposition \bproof
Since $E$ is nef, the top self intersection number $c_1^n(E)\geq 0$. If $c_1^n(E)=0$, then all degree $2n$ Chern numbers of $E$ are zero. In particular, $s_n(E)=0$. It is a contradiction by
Lemma \ref{bigbundle}. Hence the top self intersection number $c_1^n(E)>0$. Since $\det(E)$ is nef and $c_1^n(\det (E))>0$, $\det(E)$ is a big line bundle. \eproof

\bcorollary If $X$ is a compact K\"ahler manifold with nef and big tangent bundle, then $X$ is Fano. \ecorollary \bproof By Proposition \ref{bundle}, $K_X^{-1}$ is nef and big. Since $TX$ is
nef, we know $K_X^{-1}$ is ample, i.e. $X$ is Fano. \eproof

\noindent By comparing Theorem \ref{thm00} with Proposition \ref{bundle}, one may ask the following question: for an abstract vector bundle $E$, if $E$ is nef (or semi-positive ) and
$\det(E)$ is big (or ample), is $E$ big? We have a negative answer to this question.

\bexample\label{counter} On $\P^2$, let $E=T\P^2\ts \sO_{\P^2}(-1)$ be the hyperplane bundle. Then $E$ is  semi-positive in the sense of Griffiths, and $\det(E)$ is ample, but $E$ is not a
big vector bundle.\eexample\bproof By using the Hermitian metric on $E$ induced by the Fubini-Study metric, it is easy to see  that $E$ is a semi-positive vector bundle and so it is nef.
Indeed, $T\P^2$ has curvature tensor
$$R_{i\bar j k\bar \ell}=g_{i\bar j}g_{k\bar \ell}+g_{i\bar \ell} g_{k\bar j}$$
and so $E$ has curvature tensor $R^E_{i\bar j k\bar \ell}=g_{i\bar \ell}g_{k\bar j}$ where $k$ and $\ell$ are indices along the vector bundle $E$. On the other hand, $\det(E)=\sO_{\P^2}(1)$
is ample and so $c_1^2(E)=1$. However, $E$ is not a big vector bundle. Since
$$c_2(T\P^2)=c_2(E\ts
\sO_{\P^2}(1))=c_2(E)+c_1^2(\sO_{\P^2}(1))+c_1(E)\cdot c_1(\sO_{\P^2}(1))=3,$$ we have
 $c_2(E)=1,$ and so
$$s_2(E)=c_1^2(E)-c_2(E)=0.$$
Therefore,  by Lemma \ref{bigbundle}, $E$ is not a big vector bundle.\eproof

\begin{question} Suppose $X$
is a Fano manifold with nef tangent bundle. Is the (abstract) vector bundle $TX$ semi-positive in the sense of Griffiths? Is $TX$ a big vector bundle?\\

\noindent For example, $\P(T^*\P^n)$ ($n\geq 2$) is a Fano manifold with nef tangent bundle since it is homogeneous. When $n=2,$ $\P(T^*\P^2)$  has big and semi-ample tangent bundle by
Theorem \ref{SW}.
 It is also known
that $\P(T^*\P^n)$ does not admit a  smooth \textbf{K\"ahler metric} with semi-positive holomorphic bisectional curvature according to the classification in Theorem \ref{main00}. However, it
is not clear whether the tangent bundle of $\P(T^*\P^n)$ is semi-positive in the sense of Griffiths, or equivalently, whether it has a smooth \textbf{Hermitian metric} with semi-positive
holomorphic bisectional curvature. When $n>2$, is the tangent bundle of $\P(T^*\P^n)$ big?
\end{question}

 As motivated by these questions, in the next section, we
investigate compact complex manifolds with semi-positive tangent bundles.

\section{Complex manifolds with semi-positive tangent
bundles}\label{semi}

In this section, we study complex manifolds with semi-positive tangent bundles. Suppose the abstract tangent bundle $TX$ has a smooth Hermitian metric $h$ with semi-positive curvature in the
sense of Griffiths, or equivalently, $(X,h)$ is a Hermitian manifold with semi-positive holomorphic bisectional curvature.

\btheorem\label{thm11}  Let $(X,\omega)$ be a compact Hermitian manifold with semi-positive holomorphic bisectional curvature, then $\kappa(X)\leq 0$, and  either \bd\item
$\kappa(X)=-\infty$; or
\item $X$ is a complex parallelizable manifold. Moreover, $(X,\omega)$ has flat curvature and $d^*\omega=0$. \ed \etheorem

\bremark A complex manifold $X$ of complex dimension $n$ is called complex parallelizable if there exist $n$  holomorphic vector fields linearly independent everywhere.
  Note that every complex parallelizible manifold has a
balanced Hermitian metric with flat curvature tensor and the canonical line bundle is holomorphically trivial, and so the Kodaira dimension is zero. It is proved by Wang in
\cite[Corollary~2]{Wang} that a complex parallelizable manifold is K\"ahler if and only if it is a torus.\eremark

\bproof Since $(X,\omega)$ has semi-positive holomorphic bisectional curvature, $K^{-1}_X$ is  semi-positive and so nef. Suppose $\kappa(X)\geq 0$, i.e. there exists some positive integer $m$
such that $H^0(X,K^{\ts m}_X)$ has a non zero element $\sigma$. Then $\sigma$ does not vanish everywhere (\cite[Proposition~1.16]{DPS}). In particular, $K_X^{\ts m}$ is a holomorphically
trivial line bundle, i.e. $K_X^{\ts m}=\sO_X$. In this case, we obtain $\kappa(X)=0$. Let $h$ be the trivial Hermitian metric on $K_X^{\ts m}$, i.e. $\sq\p\bp\log h=0$. On the other hand,
$K_X^{\ts m}$ has a smooth Hermitian metric $\frac{1}{\left[\det(\omega)\right]^m}$. Hence, there exists a positive smooth function $\phi\in C^\infty(X)$ such that \beq
\frac{1}{\left[\det(\omega)\right]^m}=\phi\cdot h, \eeq and the line bundle $K_X^{\ts m}$ has curvature form
$$-mRic(\omega)=-\sq\p\bp\log h-\sq\p\bp\log\phi=-\sq\p\bp\log \phi\leq 0.$$
By maximum principle, we know $\phi$ is  constant. Therefore $Ric(\omega)=0$. Since $(X,\omega)$ has semi-positive holomorphic bisectional curvature, we know $R_{i\bar j k\bar \ell}=0$.
Indeed, without loss of generality, we assume $g_{i\bar j}=\delta_{ij}$ at a fixed point $p\in X$, and so the Ricci curvature has components $R_{i\bar j}=\sum_{k=1}^nR_{i\bar j k \bar k}=0$.
If we choose $b=(1,0,\cdots,0)$, then for any $a\in \C^n$, we have $R_{i\bar j k\bar \ell} a^i\bar a^j b^k\bar b^\ell=R_{i\bar j 1\bar 1}a^i\bar a^j\geq 0$. Similarly, we have $R_{i\bar j
k\bar k} a^i\bar a^j\geq 0$ for all $k=1,\cdots, n$. By the Ricci flat condition, we have $R_{i\bar j k\bar k}a^i\bar a^j=0$ for all $a\in \C$ and $k=1,\cdots, n$. We deduce $R_{i\bar j k\bar
k}=0$ for any $i,j,k$. Now for any $a\in \C^n$, we define  $H_{k\bar \ell}=R_{i\bar j k\bar \ell} a^i\bar a^j$. Then $H=(H_{k\bar \ell})$ is a  semi-positive Hermitian matrix. Since $tr H=0$,
$H$ is the zero matrix. That is, for any $a\in \C^n$ and $k,\ell$,  we have $R_{i\bar j k\bar \ell} a^i\bar a^j=0$. Finally, we obtain $R_{i\bar j k\bar \ell}=0$. Since $(X,\omega)$ is
Chern-flat,  $X$ is a complex parallelizable manifold (e.g. \cite[Proposition~2.4]{DLV} and \cite{AS}).
 On the other hand, it is well-known that if $(X,\omega)$ is Chern-flat, $d^*\omega=0$
 (e.g.,\cite[Corollary~2]{LYZ}).
 \eproof

\noindent The following application of Theorem \ref{thm11} will be used frequently. \bcorollary\label{flat} Let $(X,\omega)$ be a compact Hermitian surface. If $(X,\omega)$ has semi-positive
holomorphic bisectional curvature and $K_X$ is a holomorphic torsion, i.e. $K_X^{\ts m}=\sO_X$ for some integer $m\in \N^+$, then $(X,\omega)$ is a torus. \ecorollary \bproof Since
$\kappa(X)=0$, as shown in the proof of Theorem \ref{thm11}, $(X,\omega)$ is a parallelizable complex surface with $d^*\omega=0$. Since $\dim_\C X=2$, $d^*\omega=0$ implies $dw=0$, i.e.
$(X,\omega)$ is K\"ahler. Hence $(X,\omega)$ is a flat torus. \eproof

Now we are ready to classify compact complex surfaces with semi-positive tangent bundles. Note that, we only assume $X$ has a \textbf{Hermitian metric} with semi-positive holomorphic
bisectional curvature.

\btheorem Let $X$ be a compact K\"ahler surface. If $TX$ is (Hermitian) semi-positive, then $X$ is one of the following: \bd
\item $X$ is a torus;
\item $X$ is $\P^2$;
\item $X$ is $\P^1\times \P^1$;
\item $X$ is a ruled surface over an elliptic curve $C$, and $X$ is covered $\C\times \P^1$.
\ed \etheorem \bproof Suppose $TX$ is semi-positive. If $X$ is not a torus, then by Theorem \ref{thm11}, $\kappa(X)=-\infty$. Let $X_{\text{min}}$ be a minimal model of $X$. Since
$\kappa(X_{\text{min}})=-\infty$, by Kodaira-Enriques classification, $X_{\text{min}}$ has algebraic dimension $2$ and so $X_{\text{min}}$ is projective. Therefore, $X$ is also projective. By
\cite[Proposition~2.1]{CP}, $X$ is minimal, i.e. $X=X_{\text{min}}$ since $X$ has nef tangent bundle. By Campana-Peterell's classification of projective surfaces with nef tangent bundles
(\cite[Theorem~3.1]{CP}), $X$ is one of the following \bd
\item $X$ is an abelian surface;
\item $X$ is a hyperelliptic surface;
\item $X=\P^2$;
\item $X=\P^1\times \P^1$;
\item $X=\P(E^*)$ where $E$ is a rank $2$-vector bundle on an elliptic curve $C$ with either
\bd \item $E=\sO_C\ds L$, with $\deg(L)=0$; or
\item $E$ is given by a non-split extension $0\>\sO_C\>E\>L\>0$ with $L=\sO_C$ or $\deg L=1$.
\ed \ed It is obvious that  abelian surfaces, $\P^2$, $\P^1\times\P^1$ all have canonical Hermitian metrics with semi-positive holomorphic bisectional curvature. By Corollary \ref{flat}, a
hyperelliptic surface can not admit a Hermitian metric with semi-positive holomorphic bisectional curvature since its canonical line bundle is a torsion.  Next, we show that, in case $(5)$,
if $X=\P(E^*)$ has semi-positive tangent bundle, then $X$ is a ruled surface over an elliptic curve $C$ which is covered by $\C\times \P^1$. Indeed, by the exact sequence
$0\>T_{X/C}\>TX\>\pi^*(TC)\>0$ where $\pi:X\>C$, we obtain the dual sequence \beq 0\>\pi^*\sO_C\>T^*X\>T^*_{X/C}\>0, \eeq since $TX=\sO_C$. Suppose $TX$ is semi-positive in the sense of
Griffiths, $T^*X$ is semi-negative in the sense of Griffiths. It is well-known that, the holomorphic bisectional curvature is decreasing in subbundles, and so the induced Hermitian metric on
the subbundle $\pi^*\sO_C$ also has semi-negative curvature in the sense of Griffiths (\cite[p.~79]{GH}). Since the line bundle $\pi^*\sO_C$ is trivial, that induced metric on $\pi^*\sO_C$
must be flat by maximum principle. In particular, the second fundamental form of $\pi^*\sO_C$ in $T^*X$ is zero. Therefore,   the Hermitian metric on $TX$ splits into a direct product and the
tangent bundle $TX$ splits into the {holomorphic} direct sum
$$TX=\pi^*\sO_C\ds T_{X/C}.$$ We deduce $X$ is a ruled surface
over an elliptic curve $C$ covered by $\tilde X=\C\times \P^1$. Or equivalently, $X=\P(E^*)$ with $E=\sO_C\ds L$ where $\deg(L)=0$ on $C$. Moreover, it is also well-known that for the
non-split extension, the ant-canonical line bundle $K_X^{-1}$ of $X=\P(E^*)$ can not be semi-positive \cite[Example~3.5]{DPS}. \eproof

\noindent In the following, we classify non-K\"ahler surfaces with semi-positive tangent bundles.
 \btheorem Let $(X,\omega)$ be a non-K\"ahler compact
complex surface with semi-positive holomorphic bisectional curvature. Then $X$ is a Hopf surface. \etheorem \bproof Suppose $X$ is a non-K\"ahler complex surface. By Theorem \ref{thm11}, we
have $\kappa(X)=-\infty$ since  when $\kappa(X)=0$, $(X,\omega)$ is balanced and so it is K\"ahler. By the Enriques-Kodaira classification, the minimal model $X_{\text{min}}$ of $X$ is a
VII$_0$ surface, i.e. $X$ is obtained from $X_{\text{min}}$ by successive blowing-ups.

    We  give a  straightforward proof that if $(X,\omega)$ has semi-positive holomorphic
    bisectional curvature, then $X$ is minimal, i.e.
    $X=X_{\text{min}}$. Here we can not use methods in algebraic
    geometry since the ambient manifold is  non-K\"ahler and the curvature condition may not be preserved
    under birational maps, finite \'etale covers, blowing-ups, or blowing-downs(cf.\cite[Proposition~6.3]{DPS}).
 By definition,
$X_{\text{min}}$ is a compact complex surface with $b_1(X_{\text{min}})=1$ and $\kappa(X_\text{min})=-\infty$. It is well-known that the first Betti number $b_1$ is invariant under
blowing-ups, i.e. $b_1(X)=1$. By \cite[Theorem ~2.7 on p.139]{BHPV}, we know $b_1(X)=h^{1,0}(X)+h^{0,1}(X)$ and $h^{1,0}(X)\leq h^{0,1}(X)$, hence $h^{0,1}(X)=1$. Since $\kappa(X)=-\infty$,
we have $h^{0,2}(X)=h^{2,0}(X)=h^0(X,K_X)=0$. Therefore, by the Euler-Poincar\'e characteristic formula, we get
$$\chi(\sO_X)=1-h^{0,1}(X)+h^{0,2}(X)=0.$$ On the other hand, by the
Noether-Riemann-Roch formula,
$$\chi(\sO_X)=\frac{1}{12}(c_1^2(X)+c_2(X))=0,$$ we have $c_2(X)=-c_1^2(X)$.
$c_2(X)$ is also the Euler characteristic $e(X)$ of $X$, i.e.
$$c_2(X)=e(X)=2-2b_1(X)+b_2(X)=b_2(X)$$
and so $c_1^2(X)=-b_2(X)\leq 0$. Since $(X,\omega)$ has semi-positive holomorphic bisectional curvature, we obtain $c_1^2(X)\geq 0$. Hence $c_2(X)=b_2(X)=0$. On the other hand, blowing-ups
increase the second Betti number at least by $1$. We conclude that $X=X_{\min}$.

  Hence, $X$ is a VII$_0$ surface with $b_2(X)=0$. By
  Kodaira-Enriques's classification (see also \cite{LYZ2}), $X$ is either
 \bd \item a  Hopf surface (whose universal cover is  $\C^2\setminus\{0\}$); or

 \item  an Inoue surface, i.e. $b_1(X)=1$, $b_2(X)=0$ and
 $\kappa(X)=-\infty$, without any curve.
\ed As shown in \cite[Proposition~6.4]{DPS}, the holomorphic tangent bundles of Inoue surfaces are not nef. In particular, Inoue surfaces can not admit smooth Hermitian metrics with
semi-positive holomorphic bisectional curvature. Finally, we deduce that $X$ is a Hopf surface.
 \eproof

 \noindent A compact complex surface $X$ is called a Hopf surface if
its universal covering is analytically isomorphic to $\C^2\setminus\{0\}$. It has been prove by Kodaira that its fundamental group $\pi_1(X)$ is a finite extension of an infinite cyclic group
generated by a biholomorphic contraction which takes the form \beq (z,w)\>(az,bw+\lambda z^m) \eeq where $a,b,\lambda \in \C$, $|a|\geq |b|>1$, $m\in \N^*$ and $\lambda(a-b^m)=0$. Hence,
there are two different cases:

\bd \item the Hopf surface $H_{a,b}$ of class $1$ if $\lambda=0$;
\item the Hopf surface $H_{a,b,\lambda,m}$ of class $0$ if $\lambda\neq
0$ and $a=b^m$. \ed

\noindent In the following, we consider the Hopf surface of class $1$. Let $H_{a,b}=\C^{2}\setminus\{0\}/\sim$ where $(z,w)\sim(az,bw)$ and $|a|\geq |b|> 1$. We set $k_1=\log|a|$ and
$k_2=\log|b|$. Define a real smooth function \beq \Phi(z,w)=e^{\frac{k_1+k_2}{2\pi}\theta}\eeq where $\theta(z,w)$ is a real smooth function defined by \beq
|z|^2e^{-\frac{k_1\theta}{\pi}}+|w|^2e^{-\frac{k_2\theta}{\pi}}=1.\label{0}\eeq This is well-defined since for fixed $(z,w)$ the function $t\>|z|^2|a|^t+|w|^2|b|^t$ is strictly increasing
with image $\R_+$(\cite{GO}). Let $\alpha=\frac{2k_1}{k_1+k_2}$ and so $1\leq\alpha<2$. Then the key equation (\ref{0}) is equivalent to \beq
|z|^2\Phi^{-\alpha}+|w|^2\Phi^{\alpha-2}=1.\label{ke} \eeq It is easy to see that
$$\theta(az,bw)=\theta(z,w)+2\pi, \qtq{and}\Phi(az,bw)=|a||b|\Phi(z,w).$$
When $\alpha=1$, i.e. $|a|=|b|$, we have \beq \Phi=|z|^2+|w|^2.\eeq

\blemma\label{aa} $\displaystyle |z|^2\Phi^{-\alpha}$ and $|w|^2\Phi^{\alpha-2}$ are well-defined on $H_{a,b}$. \elemma \bproof Indeed, $$ |az|^2\Phi^{-\alpha}(az,bw)=|a|^2
|a|^{-\alpha}|b|^{-\alpha} \cdot |z|^2\Phi^{-\alpha}(z,w) $$  and
$$|a|^2 |a|^{-\alpha}|b|^{-\alpha}=e^{k_1(2-\alpha)} e^{-k_2\alpha}=1.$$
Similarly, we can show $|w|^2\Phi^{2-\alpha}$ is well-defined on $H_{a,b}$. \eproof

\noindent By Lemma \ref{aa}, we know \beq \omega=\sq \left(\lambda_1\Phi^{-\alpha}dz\wedge d\bar z+\lambda_2\Phi^{\alpha-2}dw\wedge d\bar w\right)\label{gmetric}\eeq is a well-defined
Hermitian metric on $H_{a,b}$ for any $\lambda_1,\lambda_2\in \R^+$. It is easy to see that the (first Chern) Ricci curvature of $\omega$ is \beq Ric(\omega)=-\sq\p\bp\log\det(\omega)=2\sq
\p\bp\log \Phi.\eeq The next lemma shows $Ric(\omega)\geq 0$ and $Ric(\omega)\wedge Ric(\omega)=0$. \blemma\label{bb0} $\sq\p\bp\log\Phi$ has a semi-positive matrix representation
\beq\frac{\Phi^{-2}}{\Delta^3}\left[\begin{array}{lcr}  (\alpha-2)^2|w|^2&\alpha(\alpha-2)\bar w z\\ \alpha(\alpha-2)\bar z w& \alpha^2|z|^2
\end{array}\right],\label{222}\eeq
and $\sq\p\Phi\wedge \bp\Phi$ has a matrix representation \beq \frac{1}{\Phi^{2\alpha-2}\Delta^2}\left[\begin{array}{lcr} |z|^2&\bar w z\Phi^{2\alpha-2}\\ \bar z w \Phi^{2\alpha-2}&
|w|^2\Phi^{4\alpha-4}
\end{array}\right],\label{10}\eeq
where $\Delta$ is a globally defined function on $H_{a,b}$ given by \beq \Delta=\alpha|z|^2 \Phi^{-\alpha}+(2-\alpha)|w|^2\Phi^{\alpha-2}.\label{delta0}\eeq
 In particular, $(\sq\p\bp\log\Phi)^2=0$. \elemma
\bproof It is proved in the Appendix. \eproof

\bproposition\label{hab} On every Hopf surface $H_{a,b}$, there exists a Gauduchon metric with semi-positive holomorphic bisectional curvature. \eproposition

\bproof We show that  \beq \omega=\sq \left(\frac{\Phi^{-\alpha} }{\alpha^2}dz\wedge d\bar z+\frac{\Phi^{\alpha-2}}{(2-\alpha)^2}dw\wedge d\bar w\right)\eeq is
 a Gauduchon metric with semi-positive holomorphic bisectional
 curvature.

   At first, we show $\omega$ is Gauduchon, i.e.  $\p\bp\omega=0$. Indeed, by the elementary
identity $\p\bp f=f\p\bp\log f+f^{-1}\p f\wedge \bp f$, we obtain
$$\p\bp\Phi^\mu=\mu\Phi^\mu\p\bp\log\Phi+\mu^2\Phi^{\mu-2}\p\Phi\wedge \bp\Phi.$$
In particular we have \be \p_w\p_{\bar w}\Phi^{-\alpha}&=&-\alpha\Phi^{-\alpha}\p_w\p_{\bar w}\log\Phi+\alpha^2\Phi^{-\alpha-2}\p_w\Phi\cdot \p_{\bar w}\Phi\\ &=&-\alpha\Phi^{-\alpha}\cdot
\frac{\Phi^{-2}}{\Delta^3} (\alpha^2 |z|^2)+\alpha^2\Phi^{-\alpha-2}\cdot \frac{|w|^2\Phi^{4\alpha-4}}{\Phi^{2\alpha-2}\Delta^2}\ee where we use (\ref{222}) and (\ref{10}) in the second
identity.  Hence \beq \p_w\p_{\bar w}\left(\frac{\Phi^{-\alpha}}{\alpha^2}\right)=\frac{-\alpha|z|^2\Phi^{-2-\alpha}}{\Delta^3}+\frac{|w|^2\Phi^{\alpha-4}}{\Delta^2}.\eeq Similarly, we have
\beq \p_z\p_{\bar z}\left(\frac{\Phi^{\alpha-2}}{(\alpha-2)^2}\right)=-\frac{(2-\alpha)|w|^2\Phi^{\alpha-4}}{\Delta^3}+\frac{|z|^2\Phi^{-\alpha-2}}{\Delta^2}.\eeq Now it is obvious that
$$\p_w\p_{\bar
w}\left(\frac{\Phi^{-\alpha}}{\alpha^2}\right)+ \p_z\p_{\bar z}\left(\frac{\Phi^{\alpha-2}}{(\alpha-2)^2}\right)=0$$ where we use equations (\ref{ke}) and $(\ref{delta0})$. This implies
$\p\bp\omega=0$.

Next, we prove $\omega$ has semi-positive holomorphic bisectional curvature. To write down the holomorphic bisectional curvature, we introduce new notations, $z^1=z$ and $z^2=w$. Moreover,
let $\omega=\sq \sum_{i,j=1}^2g_{i\bar j}dz^i\wedge d\bar z^j$ with $g_{ij}=f_i \delta_{ij}$, where $f_1=\frac{\Phi^{-\alpha}}{\alpha^2}$ and $ f_2=\frac{\Phi^{\alpha-2}}{(\alpha-2)^2}$ or
equivalently \beq f_i=\frac{\Phi^{(2i-3)\alpha+2(1-i)}}{((2i-3)\alpha+2(1-i))^2},\ \ \ i=1, 2.\eeq Therefore, the Christoffel  symbols of $\omega$ are \be
\Gamma_{ik}^{p}&=&\sum_{q=1}^2g^{p\bar q}\frac{\p g_{k\bar q}}{\p z^i}=\frac{\p\log f_k}{\p z^i}\delta_{kp}=\frac{\p\log\Phi}{\p z^i}\cdot ((2k-3)\alpha+2(1-k))\delta_{kp}.\ee Hence $
R_{i\bar jk}^q=-\p_{\bar j}\Gamma_{ik}^p=\frac{\p^2\log\Phi}{\p z^i\p\bar z^j}\cdot ((3-2k)\alpha+2(k-1))\delta_{kp}, $ and \beq R_{i\bar j k\bar\ell}=\frac{\p^2\log\Phi}{\p z^i\p\bar
z^j}\cdot ((3-2k)\alpha+2(k-1)) f_k\delta_{k\ell}.\eeq Note that, by Lemma \ref{bb0}, $\left(\frac{\p^2\log\Phi}{\p z^i\p\bar z^j}\right)$ is semi-positive and \be\left(((3-2k)\alpha+2(k-1))
f_k\delta_{k\ell}\right)=\left[\begin{array}{cc}\frac{\Phi^{-\alpha}}{\alpha}&0\\ 0&\frac{\Phi^{\alpha-2}}{2-\alpha}
\end{array}\right].
\ee We deduce $R_{i\bar j k\bar\ell}$ is semi-positive in the sense of Griffiths, i.e. $\omega$ has semi-positive holomorphic bisectional curvature. \eproof

Let $X$ be a  complex manifold. $X$ is said to be a complex Calabi-Yau manifold if $c_1(X)=0$.

\bcorollary\label{CY} Let $X$ be  a compact complex Calabi-Yau manifold in the Fujiki class $\mathscr C$ (class of manifolds bimeromorphic to K\"ahler manifolds). Suppose $X$ has a Hermitian
metric $\omega$ with semi-positive holomorphic bisectional curvature, then $X$ is a torus. \ecorollary

\bproof Let $X$ be a compact Calabi-Yau manifold in the class $\mathscr C$, then by a result of \cite[Theorem~1.5]{T}, $K_X$ is a holomorphic torsion, i.e. there exists a positive integer $m$
such that $K_X^{\ts m}=\sO_X$. Suppose $X$ has a smooth Hermitian metric $\omega$ with semi-positive holomorphic bisectional curvature, then by Theorem \ref{thm11}, $X$ is a complex
parallelizable manifold. On the other hand, by
 \cite[Corollary~1.6]{Dem91} or \cite[Proposition~3.6]{DPS}, $X$ is
 K\"ahler since $X$ is in the Fujiki class $\mathscr C$ and $TX$ is nef. It is well-known that a complex parallelizable
 manifold is  K\"ahler if and only if it
 is a torus.
\eproof

\bremark As shown in Proposition \ref{hab}, the Hopf surface $H_{a,b}$ (and every diagonal Hopf manifold \cite{LY14}) has a Hermitian metric with semi-positive holomorphic bisectional
curvature. Since $b_2(H_{a,b})=b_2(\S^1\times \S^3)=0$, we see $c_1(H_{a,b})=0$ and so $H_{a,b}$ is a non-K\"ahler Calabi-Yau manifold. Hence, the Fujiki class condition  in Corollary
\ref{CY} is necessary. \eremark

To end this section, we give  new examples on K\"ahler and non-K\"ahler manifolds whose tangent bundles or anti-canonical line bundles are \emph{nef but not semi-positive}.

\bcorollary Let $X$ be a Kodaira surface or a hyperelliptic surface.

 \bd\item The  tangent bundle $TX$  is nef
but not semi-positive (in the sense of Griffiths);

\item The anti-canonical line
bundle of the projective bundle $\P(T^*X)$ is nef,  but it is neither semi-positive nor big.

 \ed\ecorollary
\bproof Suppose $X$ is a Kodaira surface.\\
 $(1).$ By the fibration structure $0\>T_{X/C}\>TX\>\pi^*TC\>0$ of a
 Kodaira surface, we know $TX$ is nef.
  Since the canonical line bundle of every
Kodaira surface is a torsion, i.e. $K_X^{\ts m}=\sO_X$ with $m=1,2,3,4$ or $6$, by Corollary \ref{flat}, $TX$ can not be semi-positive.

 For part $(2)$,
 let $Y:=\P(T^*X)$ and $\sO_Y(1)$ be the
tautological line bundle of $Y$ and $\pi:Y\>X$ the canonical projection. Suppose $TY$ is big, then $K_Y^{-1}=\sO_Y(2)$ is also a big line bundle. Therefore $Y$ is a Moishezon manifold with
nef tangent bundle, and so $Y$ is projective. On the projective manifold $Y$, $K_Y^{-1}$ is nef and big, and so by Kawamata-Reid-Shokurov's base point free theorem, $K_Y^{-1}$ is semi-ample.
Moreover, since $K_Y^{-1}$ is big,  $\int_Y c_1^3(Y)>0$. It implies $K_Y^{-1}$ is ample. Therefore, $\sO_Y(1)$ is ample and so is $TX$ which is a contradiction.

 Let $E=TX$. Then $\det E=K_X^{-1}$ is a holomorphic
torsion. By Corollary \ref{directimage2}, $E$ is semi-positive in the sense of Griffiths if and only if $\sO_Y(1)$ is semi-positive. Since $K_Y^{-1}=\sO_Y(2)$, and $E=TX$ can not be
semi-positive, we deduce $K_Y^{-1}$ can not be semi-positive.

 When $X$ is a hyperelliptic surface, the proof is
similar.
 \eproof

\bremark It is not clear where $\P(T^*\P^2)$ has a Hermitian metric with semi-positive holomorphic bisectional. Note that the tangent bundle of $\P(T^*\P^2)$ is semi-ample. It is related  to
a weak version of Griffiths' conjecture: if $E$ is semi-ample, then $E$ has a Hermitian metric with semi-positive curvature in the sense of Griffiths.  On the other hand, it is known that
$E\ts \det E$ has a metric with semi-positive curvature, and for large $k$, $S^kE$ has a Hermitian metric with Griffiths semi-positive curvature.

 \eremark

\section{Projective bundle $\P(T^*X)$ with nef tangent bundle
}\label{PTX} In this section, we study complex manifolds of the form $\P(T^*X)$ which also have nef tangent bundles. At first, we introduce the (maximum) irregularity of a compact complex
manifold $M$, \beq \tilde q(M)=\sup \{q(\tilde M)\ |
 \exists\ \text{a finite \'etale cover $f: \tilde M \> M$}\}, \eeq where $q(N)=h^1(N,\sO_N)$ for any complex
manifold $N$.

It is well-known that $\P(T^*\P^n)$ is homogeneous, and its tangent bundle is nef.  We have a similar converse statement and yield another characterization of $\P^n$.
 \bproposition Let $X$ be a Fano manifold of complex
dimension $n$. Suppose $\P(T^*X)$ has nef tangent bundle, then $X\cong \P^n$. \eproposition \bproof Let $Y=\P(T^*X)$ and $\pi:Y\>X$ be the projection. It is obvious that $\pi$ has fiber
$F=\P^{n-1}$. Since $F$ and $X$ are Fano manifolds,  $\tilde q(X)=q(X)=0$ and $\tilde q(F)=q(F)=0$. Therefore, from the relation (\cite[Proposition~3.12]{DPS}) \beq \tilde q(Y)\leq \tilde
q(X)+\tilde q(F),\eeq we obtain $\tilde q(Y)=0$. We claim $Y$ is Fano. Indeed, since $TY$ is nef, $c^{2n-1}_1(Y)\geq 0$. Suppose $c_1^{2n-1}(Y)=0$, then by \cite[Proposition~3.10]{DPS}, there
exists a finite \'etale cover $\tilde Y$ of $Y$ such that $q(\tilde Y)>0$ which is a contradiction since $\tilde q(Y)=0$. Hence, we have $c_1^{2n-1}(Y)>0$, i.e. $K_Y^{-1}$ is nef and big. Now
we deduce $Y$ is projective and $K_Y^{-1}$ is ample. By the adjunction formula, $K_Y^{-1}=\sO_Y(2)$. We obtain $\sO_Y(1)$ and so $TX$ are ample. Hence $X=\P^n$ by Mori's result. \eproof

In the rest of this section, we classify complex $3$-folds of the form $\P(T^*X)$ whose tangent bundles are nef.

\bproposition\label{surface} Let $X$ be a compact K\"ahler surface. If the projective bundle $\P(T^*X)$ has nef tangent bundle, then $X$ is exactly one of the following: \bd \item $X\cong
\T^2$, a flat torus;
\item $X\cong \P^2$;
\item $X$ is a hyperelliptic surface;
 \ed \eproposition

\bproof  Let $Y=\P(E^*)$ and $\pi:Y\>X$ the canonical projection. Consider the exact sequence
$$0\>T_{Y/X}\>TY\>\pi^*TX\>0.$$
Since, $TY$ is nef, the quotient bundle $\pi^*TX$ is nef (\cite[Proposition~1.15]{DPS}). On the other hand, since $\pi:Y\>X$ is a surjective holomorphic map with equidimensional fibers, we
deduce $TX$ is nef. Then $X$ is one of the following \bd
\item $X$ is a torus;
\item $X$ is a hyperelliptic surface;
\item $X=\P^2$;
\item $X=\P^1\times \P^1$
\item $X=\P(E^*)$ where $E$ is a rank $2$-vector bundle on an elliptic curve $C$ with either
\bd \item $E=\sO_C\ds L$, with $\deg(L)=0$; or
\item $E$ is given by a non-split extension $0\>\sO_C\>E\>L\>0$ with $L=\sO_C$ or $\deg L=1$.
\ed \ed

\noindent It is obvious that torus and $\P^2$ satisfy the requirement. By Corollary \ref{pn}, we can rule out $X=\P^1\times \P^1$ since $TY$ can not be nef. Now we verify that when $X$ is a
hyperelliptic surface, both $TX$ and $\P(T^*X)$ have nef tangent bundles. It is well-known that every hyperelliptic surface $X$ is a projective manifold, which admits a locally trivial
fibration $\pi:X\>C$ over an elliptic curve $C$, with an elliptic curve as a typical fiber. Moreover, $K_X$ is a torsion line bundle \cite[p.~245]{BHPV}, i.e. $K_X^{\ts m}=0$ for $m=2,3,4$,
or $6$. By the exact sequence $0\>T_{X/C}\>TX\>\pi^*TC\>0,$ we know $TX$ is nef since both $\pi^*TC$ and $T_{X/C}=K^{-1}_X\ts K_C$ are nef line bundles. Let $Y=\P(T^*X)$ and $\pi_1:Y\>X$.
Then $Y$ is a $\P^1$-bundle over $X$. Similarly, from the exact sequence $0\>T_{Y/X}\>TY\>\pi_1^*TX\>0$ we can also deduce $TY$ is nef. Here, we only need to show $T_{Y/X}$ is nef. Indeed,
$K_Y^{-1}=\sO_Y(2)$ where $\sO_Y(1)$ is the tautological line bundle of $Y=\P(T^*X)$. Since $TX$ is nef, we know $\sO_Y(1)$ and so $K_Y^{-1}$ are nef. Since $K_X$ is a torsion line bundle and
$T_{Y/X}=K^{-1}_Y\ts \pi_1^*(K_X)$, we conclude $T_{Y/X}$ is a nef line bundle.
 If $X=\P(E^*)$ in $(5)$, then we know $Y=\P(T^*X)\>X\>C$ is a $\P^1\times
 \P^1$ bundle over $C$ since $TY$ is nef (\cite[Lemma~9.3]{CP}).  It
 is easy to see that the fiber of $Y\>C$ is isomorphic to the second Hirzrbruch surface
 $\P(\sO_{\P^1}\ds\sO_{\P^1}(-2))$. Indeed, for any $s\in C$, the
 fiber $X_s$ of $X\>C$ is $\P^1$. From the exact sequence
 $0\>T\P^1\>TX|_{\P^1}\>N_{\P^1/X}=\sO_{\P^1}\>0$, we see $TX|_{\P^1}=\sO_{\P^1}\ds
 T\P^1$. Hence, the fiber $Y_s$ of $Y\>C$ is isomorphic to $\P(T^*Y|_{\P^1})\cong
 \P(\sO_{\P^1}\ds\sO_{\P^1}(-2))$. Suppose $Y$ has nef tangent bundle, so is the fiber $
 \P(\sO_{\P^1}\ds\sO_{\P^1}(-2))$(\cite[Proposition~2.1]{CP}). However, the second Hirzebruch
 surface contains a $(-2)$-curve, the tangent bundle can not be nef.
 \eproof

\bproposition Let $X$ be a  non-K\"ahler  compact complex surface. If the projective bundle $\P(T^*X)$ has nef tangent bundle, then either \bd

\item $X$ is a Kodaira surface; or

\item $X$ is a Hopf surface.

\ed\eproposition

\bproof By similar arguments as in the proof of Proposition \ref{surface}, we deduce $X$ has nef tangent bundle. It is well-known that, either \bd

\item $X$ is a Kodaira surface; or

\item $X$ is a Hopf surface.
\ed Now we verify  $\P(T^*X)$ has nef tangent bundle in both cases. Let $Y=\P(T^*X)$ and $\pi:Y\>X$. Let $\sO_Y(1)$ be the tautological line bundle of $Y$ and $\pi:Y\>X$ the canonical
projection, then by adjunction formula (\ref{adjunction}), we have $K_Y^{-1}=\sO_Y(2)$. Since $TX$ is nef, by definition, $\sO_Y(1)$ and $K_Y^{-1}$ are nef. Moreover, we have the exact
sequence $0\>T_{Y/X}\>TY\>\pi^*TX\>0,$ where $T_{Y/X}=K^{-1}_{Y/X}=\sO_{Y}(2)\ts \pi^*K_X$. To obtain the nefness of $TY$, we only need to show $\sO_{Y}(2)\ts \pi^*K_X$ is nef.

Suppose $X$ is a Kodaira surface. It is well-known that $K_X$ is a torsion, hence $\sO_{Y}(2)\ts \pi^*K_X$ is nef.

  Let $X$ be a Hopf surface. Although $c_1(K_X)=0$, $K_X$ is not a torsion. We will construct explicit Hermitian metrics on $T_{Y/X}=\sO_Y(2)\ts \pi^*(K_X)$ to
show it is a nef line bundle. As a model case, we show $T_{Y/X}$ is nef for the diagonal Hopf surface. Let $\omega=\frac{\sq(dz\wedge d\bar z+dw\wedge d\bar w)}{|z|^2+|w|^2}$
 be the standard Hermitian metric on $X$. Let $[W_1,W_2]$ be the homogeneous coordinates on the fiber of $T^*X$, then
by using the curvature formula (\ref{inducedcurvature}), the tautological line bundle $\sO_Y(1)$ has curvature \beq\sq \p\bp\log
\left((|z|^2+|w|^2)|W|^2\right)=\frac{1}{2}Ric(\omega)+\sq\p\bp\log|W|^2,\label{cur} \eeq since $Ric(\omega)=2\sq\p\bp\log(|z|^2+|w|^2)$.
 The induced metric on $T_{Y/X}=\sO_{Y}(2)\ts \pi^*K_X$ has
curvature \beq 2\left(\sq\p\bp\log |z|^2+\sq\p\bp\log|W|^2\right)-Ric(\omega)=2\sq \p\bp\log|W|^2\eeq which is the Ricci curvature of the fiber $\P^1$. Hence, $T_{Y/X}$ is semi-positive and
so nef over $Y$.

Next, on a general Hopf surface $X=H_{a,b}$( $a\neq b$), we choose a Hermitian metric on $X$ as in (\ref{gmetric})$$\omega=\sq \left(\lambda_1\Phi^{-\alpha}dz\wedge d\bar
z+\lambda_2\Phi^{\alpha-2}dw\wedge d\bar w\right).$$ Then $T_{Y/X}=\sO_{Y}(2)\ts \pi^*K_X$ has an induced metric \be &&2\sq\p\bp\log(\lambda_1^{-1}\Phi^\alpha
|W_1|^2+\lambda_2^{-1}\Phi^{2-\alpha}|W_2|^2)-2\sq\p\bp\log\Phi\\&=&2\sq\p\bp\log(\lambda_1^{-1}\Phi^{\alpha-1} |W_1|^2+\lambda_2^{-1}\Phi^{1-\alpha}|W_2|^2). \ee Fix a Hermitian metric
$\omega_Y$ on $Y$. Note that $\sq \p\bp\log(\lambda_1^{-1}\Phi^{\alpha-1} |W_1|^2)$ is semi-positive by Lemma \ref{bb0}. Hence, for any $\eps>0$, we can fix $\lambda_1$ and choose $\lambda_2$
large enough such that
$$2\sq\p\bp\log(\lambda_1^{-1}\Phi^{\alpha-1}
|W_1|^2+\lambda_2^{-1}\Phi^{1-\alpha}|W_2|^2)\geq -\eps\omega_Y.$$ For a Hopf surface of type $0$, since the $z$-direction is still invariant, we can use similar arguments as above to show
$T_{Y/X}=\sO_{Y}(2)\ts \pi^*K_X$ is nef (see also the arguments in \cite[Proposition~6.3]{DPS})
 \eproof

\section{Appendix}\label{appendix}

In this appendix, we prove Lemma \ref{bb0}, i.e.

\blemma\label{bb}$\sq\p\bp\log\Phi$ has a matrix representation \beq\frac{\Phi^{-2}}{\Delta^3}\left[\begin{array}{lcr}  (\alpha-2)^2|w|^2&\alpha(\alpha-2)\bar w z\\ \alpha(\alpha-2)\bar z w&
\alpha^2|z|^2
\end{array}\right],\label{22}\eeq
and $\sq\p\Phi\wedge \bp\Phi$ has a matrix representation \beq \frac{1}{\Phi^{2\alpha-2}\Delta^2}\left[\begin{array}{lcr} |z|^2&\bar w z\Phi^{2\alpha-2}\\ \bar z w \Phi^{2\alpha-2}&
|w|^2\Phi^{4\alpha-4}
\end{array}\right],\eeq
where $\Delta$ is a globally defined function on $H_{a,b}$ given by \beq \Delta=\alpha|z|^2 \Phi^{-\alpha}+(2-\alpha)|w|^2\Phi^{\alpha-2}.\label{delta00}\eeq
 In particular, $(\sq\p\bp\log\Phi)^2=0$. \elemma
\bproof By taking $\p$ on equation (\ref{ke}), i.e. $|z|^2\Phi^{-\alpha}+|w|^2\Phi^{\alpha-2}=1$, we obtain
$$\p|z|^2\cdot \Phi^{-\alpha}-\alpha|z|^2\Phi^{-\alpha-1}\cdot\p\Phi+\p|w|^2\cdot \Phi^{\alpha-2}+(\alpha-2)|w|^2\Phi^{\alpha-3}\cdot\p\Phi=0$$
and so \beq \p\Phi= \frac{\p|z|^2\cdot\Phi^{-\alpha}+\p|w|^2\cdot |\Phi|^{\alpha-2}}{\alpha|z|^2\Phi^{-\alpha-1}+(2-\alpha)|w|^2\Phi^{\alpha-3}}= \frac{\p|z|^2+\p|w|^2\cdot
|\Phi|^{2\alpha-2}}{\Phi^{\alpha-1}\Delta}. \eeq Similarly, we have \beq \bp\Phi= \frac{\bp|z|^2\cdot\Phi^{-\alpha}+\bp|w|^2\cdot
|\Phi|^{\alpha-2}}{\alpha|z|^2\Phi^{-\alpha-1}+(2-\alpha)|w|^2\Phi^{\alpha-3}}= \frac{\bp|z|^2+\bp|w|^2\cdot |\Phi|^{2\alpha-2}}{\Phi^{\alpha-1}\Delta}. \eeq  Their wedge product is
$$\p\Phi\wedge \bp\Phi=\frac{\p|z|^2\cdot \bp|z|^2+\p|w|^2\cdot \bp|z|^2\cdot \Phi^{2\alpha-2}+\p|z|^2\cdot
\bp|w|^2\cdot\Phi^{2\alpha-2}+\p|w|^2\cdot\bp|w|^2\cdot\Phi^{4\alpha-4}}{\Phi^{2\alpha-2}\Delta^2}$$
and in the matrix form it is \beq\p\Phi\wedge \bp\Phi\sim \frac{1}{\Phi^{2\alpha-2}\Delta^2}\left[\begin{array}{lcr} |z|^2&\bar w z\Phi^{2\alpha-2}\\ \bar z w \Phi^{2\alpha-2}&
|w|^2\Phi^{4\alpha-4}
\end{array}\right].\label{1}\eeq
Since $\bp\left(|z|^2\Phi^{-\alpha}+|w|^2\Phi^{\alpha-2}\right)=0$, i.e.
$$0=\left(|z|^2 (-\alpha)\Phi^{-1}+|w|^2(\alpha-2)\Phi^{2\alpha-3}\right)\bp\Phi+\left(\bp|z|^2+\bp|w|^2\cdot\Phi^{2\alpha-2}\right)$$
by taking $\p$ again, we have \be 0&=&\left(|z|^2 (-\alpha)\Phi^{-1}+|w|^2(\alpha-2)\Phi^{2\alpha-3}\right)\p\bp\Phi\\ &+&\left(\p|z|^2\cdot(-\alpha)\cdot\Phi^{-1}+\p|w|^2\cdot
(\alpha-2)\Phi^{2\alpha-3}\right)\wedge \bp\Phi\\&+&\left(\alpha|z|^2\Phi^{-2}+(\alpha-2)(2\alpha-3)|w|^2\Phi^{2\alpha-4}\right)\p\Phi\wedge \bp\Phi\\ &+&\p\bp|z|^2+\p\bp|w|^2\cdot
\Phi^{2\alpha-2}+(2\alpha-2)\Phi^{2\alpha-3}\p\Phi\wedge \bp|w|^2. \ee Hence, we find
 \be \p\bp\Phi&=&\frac{\p\bp|z|^2+\p\bp|w|^2\cdot
\Phi^{2\alpha-2}}{\Phi^{\alpha-1}\Delta}+\frac{\p|z|^2\cdot(-\alpha)\cdot\Phi^{-1}+\p|w|^2\cdot (\alpha-2)\Phi^{2\alpha-3}}{\Phi^{\alpha-1}\Delta}\wedge
\bp\Phi\\&&+\frac{\alpha|z|^2\Phi^{-2}+(\alpha-2)(2\alpha-3)|w|^2\Phi^{2\alpha-4}}{\Phi^{\alpha-1}\Delta}\p\Phi\wedge \bp\Phi+\frac{(2\alpha-2)\Phi^{2\alpha-3}\p\Phi\wedge
\bp|w|^2}{\Phi^{\alpha-1}\Delta}\\ &=&A+B+C+D,\ee where $A,B,C$ and $D$ are four summands in the previous line  respectively. We can simplify $A$ and write it as
$$A=\frac{( \alpha|z|^2\Phi^{-1}+(2-\alpha)|w|^2\Phi^{2\alpha-3})\p\bp|z|^2+\left( \alpha|z|^2\Phi^{2\alpha-3}+(2-\alpha)|w|^2\Phi^{4\alpha-5}\right)\p\bp|w|^2}{\Phi^{2\alpha-2}\Delta^2}$$
and the corresponding matrix form is \beq A\sim \frac{1}{\Phi^{2\alpha-2}\Delta^2}\left[\begin{array}{lcr} \alpha|z|^2\Phi^{-1}+(2-\alpha)|w|^2\Phi^{2\alpha-3}&0\\ 0&
\alpha|z|^2\Phi^{2\alpha-3}+(2-\alpha)|w|^2\Phi^{4\alpha-5}
\end{array}\right].\eeq

\noindent Similarly, $B$ has the matrix form \beq B\sim { \frac{1}{\Phi^{2\alpha-2}\Delta^2}\left[\begin{array}{ll} -\alpha|z|^2\Phi^{-1}&(\alpha-2)\bar w z \Phi^{2\alpha-3}\\ -\alpha\bar z w
\Phi^{2\alpha-3}& (\alpha-2)|w|^2\Phi^{4\alpha-5}
\end{array}\right].}\eeq
The matrix form of $C$ is \beq C\sim \frac{\alpha|z|^2\Phi^{-2}+(\alpha-2)(2\alpha-3)|w|^2\Phi^{2\alpha-4}}{\Phi^{\alpha-1}\Delta}\cdot
\frac{1}{\Phi^{2\alpha-2}\Delta^2}\left[\begin{array}{ll} |z|^2&\bar w z\Phi^{2\alpha-2}\\ \bar z w \Phi^{2\alpha-2}& |w|^2\Phi^{4\alpha-4}
\end{array}\right]. \eeq
We also have
$$D=\frac{(2\alpha-2)(\p|w|^2\cdot \bp|w|^2)\Phi^{4\alpha-5}+(2\alpha-2)\p|z|^2\cdot \bp|w|^2\cdot \Phi^{2\alpha-3}}{\Phi^{2\alpha-2}\Delta^2}$$
and its matrix form \beq D\sim  \frac{1}{\Phi^{2\alpha-2}\Delta^2}\left[\begin{array}{ll} 0&0\\ (2\alpha-2)\bar z w \Phi^{2\alpha-3}& (2\alpha-2)|w|^2\Phi^{4\alpha-5}
\end{array}\right].\eeq
It is easy to see that \beq A+B+ D= \frac{1}{\Phi^{2\alpha-2}\Delta^2}\left[\begin{array}{ll} (2-\alpha)|w|^2\Phi^{2\alpha-3}&(\alpha-2)\bar w z\Phi^{2\alpha-3}\\ (\alpha-2)\bar z w
\Phi^{2\alpha-3}& \alpha\Phi^{3\alpha-3} (2\alpha-2)|w|^2\Phi^{4\alpha-5}
\end{array}\right]\label{ABD}\eeq

\noindent
 We have
$\p\bp\log\Phi=\Phi^{-1}\p\bp\Phi-\Phi^{-2}\p\Phi\wedge \bp\Phi$ and so
$$\p\bp\log\Phi=\Phi^{-1}(A+B+D)+\Phi^{-1}(C-\Phi^{-1}\p\Phi\wedge \bp\Phi).$$
Here  the computation of $C-\Phi^{-1}\p\Phi\wedge \bp\Phi$ is a little bit easier and \be C-\Phi^{-1}\p\Phi\wedge
\bp\Phi&=&\frac{\alpha|z|^2\Phi^{-2}+(\alpha-2)(2\alpha-3)|w|^2\Phi^{2\alpha-4}-\Phi^{\alpha-2}\Delta}{\Phi^{\alpha-1}\Delta}\p\Phi\wedge \bp\Phi\\
&=&\frac{(\alpha-2)(2\alpha-2)|w|^2\Phi^{2\alpha-4}}{\Phi^{\alpha-1}\Delta}\p\Phi\wedge \bp\Phi\\
&=&\frac{(\alpha-2)(2\alpha-2)|w|^2\Phi^{2\alpha-4}}{\Phi^{\alpha-1}\Delta}\cdot\frac{1}{\Phi^{2\alpha-2}\Delta^2}\left[\begin{array}{ll} |z|^2&\bar w z\Phi^{2\alpha-2}\\ \bar z w
\Phi^{2\alpha-2}& |w|^2\Phi^{4\alpha-4}
\end{array}\right].\ee
Now by using (\ref{ABD}), we obtain (\ref{22}).
 \eproof

\end{document}